\documentclass[11pt,oneside]{amsart}
\newif\ifmarkup
\markupfalse
\ifdefined\cleanversion\markupfalse\fi

\usepackage[T1]{fontenc}
\usepackage[utf8]{inputenc}
\IfFileExists{lmodern.sty}{\usepackage{lmodern}\usepackage{microtype}}{\usepackage[expansion=false]{microtype}}
\usepackage{amsmath,amssymb,amsthm,mathtools,mathrsfs}
\usepackage{enumitem}
\usepackage{xcolor}
\usepackage{comment}
\ifmarkup
  \usepackage[a4paper,left=20mm,right=60mm,top=28mm,bottom=28mm,
              marginparwidth=49mm,marginparsep=4mm,headheight=14pt]{geometry}
\else
  \usepackage[a4paper,hmargin=33mm,vmargin=32mm,headheight=14pt]{geometry}
\fi
\usepackage[colorlinks=true,linkcolor=blue!50!black,citecolor=blue!50!black,urlcolor=blue!50!black]{hyperref}

\definecolor{rev}{RGB}{200,0,0}
\ifmarkup

  \makeatletter
  \newcommand{\del}[1]{\leavevmode{\color{rev}\@delwords#1 \@nil}}
  \def\@delwords#1 #2\@nil{\@strikebox{#1}\ifx\relax#2\relax\else
    \leaders\hrule height 0.62ex depth -0.5ex\hskip 0.33em plus 0.15em minus 0.1em\relax
    \@delwords#2\@nil\fi}
  \def\@strikebox#1{\setbox0\hbox{#1}\rlap{\rule[0.5ex]{\wd0}{0.12ex}}\box0}
  \makeatother
  \newcommand{\mnote}[1]{\marginpar{\raggedright\scriptsize\sffamily\color{rev}#1}}
  \newenvironment{added}{\color{rev}}{}
  \newcommand{\delnote}[1]{\par\noindent{\color{rev}\small[\textit{Deleted:} #1]}\par}
\else
  
  \newcommand{\del}[1]{}
  \newcommand{\mnote}[1]{}
  
  \newcommand{\delnote}[1]{}
\fi
\definecolor{newrev}{RGB}{0,70,230}
\ifmarkup
  \makeatletter
  \newcommand{\delb}[1]{\leavevmode{\color{newrev}\@delwords#1 \@nil}}
  \makeatother
  
  \newcommand{\mnoteb}[1]{\marginpar{\raggedright\scriptsize\sffamily\color{newrev}#1}}
  \newenvironment{addedb}{\color{newrev}}{}
  \newcommand{\delnoteb}[1]{\par\noindent{\color{newrev}\small[\textit{Deleted:} #1]}\par}
  \hypersetup{linkcolor=green!35!black,citecolor=green!35!black,urlcolor=green!35!black}
\else
  \newcommand{\delb}[1]{}
  
  \newcommand{\mnoteb}[1]{}
  
  \newcommand{\delnoteb}[1]{}
\fi

\newtheorem{mainthm}{Theorem}

\newtheorem{theorem}{Theorem}[section]
\newtheorem{proposition}[theorem]{Proposition}
\newtheorem{lemma}[theorem]{Lemma}
\newtheorem{corollary}[theorem]{Corollary}
\newtheorem{question}[theorem]{Question}
\newtheorem{problem}[theorem]{Problem}

\theoremstyle{definition}
\newtheorem{definition}[theorem]{Definition}

\theoremstyle{remark}
\newtheorem{remark}[theorem]{Remark}
\numberwithin{equation}{section}

\newcommand{\M}{\mathcal M}
\newcommand{\Aff}{\operatorname{Aff}}

\newcommand{\id}{\mathrm I}
\newcommand{\ty}{\tau_{\mathrm Y}}
\newcommand{\Dom}{\operatorname{Dom}}
\newcommand{\Ran}{\operatorname{Ran}}
\newcommand{\Z}{\mathcal Z}
\newcommand{\splus}{\mathbin{\hat{+}}}

\newcommand{\sdot}{\mathbin{\hat{\cdot}}}

\begin{document}

\title[On the extension of traces to affiliated operators]{On the extension of non-faithful and non-normal traces of von Neumann algebras to their affiliated operators}

\author{
\large Comlan Edmond Koudjinan\\[4pt]
\small  University of Toronto, Ontario, Canada\\[4pt]
\small \href{mailto:Koudjinanedmond@gmail.com}
{\texttt{Koudjinanedmond@gmail.com}}
}
\date{\today}

\begin{abstract}
{Let $\tau$ be a positive tracial functional on a finite von Neumann algebra $\M$, not necessarily normal. For a positive self-adjoint operator $A$ affiliated with $\M$ we define $\ty(A)=\lim_{h\downarrow0}\tau(A(\id+hA)^{-1})\in[0,\infty]$. We show that $\ty(A)=\sup_R\tau(\min(A,R\id))$. 
 We prove that $\ty$ is additive on positive affiliated operators, by a form--resolvent comparison valid for arbitrary positive traces. As a consequence, $\tau$ has a positive linear extension to the span $\mathscr D_\tau$ of the positive affiliated operators of finite extended trace, and $\ty(T^*T)$ defines a pre-Hilbert structure for every positive trace. {For every trace, $\ty(|S+T|)^{2/3}\le\ty(|S|)^{2/3}+\ty(|T|)^{2/3}$; hence the operators $T$ with $\ty(|T|)<\infty$ form a $\M$-bimodule $\mathscr L_Y^1$, on which the linear extension is defined and cyclic. It turns out that $\mathscr L_Y^1=\mathscr D_\tau$;} $d(S,T)=\ty(|S-T|)^{2/3}$ is a translation-invariant pseudometric on $\mathscr L^1_Y$ and is a metric iff $\tau$ is faithful. On the approximation side, a Chernoff family in $\M$ forces affiliation of its generator. Trace defects of Chernoff families need not converge to $\ty(A)$, even for faithful normal traces, while two-sided resolvent comparison suffices. {For faithful normal traces, $\ty$ is the canonical extension and the associated space $\mathscr L_Y^p$ coincide with the classical $L^p(\M,\tau)$ as well as their norms ($1\leq p<\infty$).}}
\end{abstract}

\begingroup
\renewcommand{\uppercasenonmath}[1]{}
\let\MakeUppercase\relax
\maketitle
\endgroup


\section{Introduction}\label{sec:intro}

\subsection{Motivation}

\smallskip\noindent\textit{Traces of unbounded operators.}
A trace $\tau$ is defined on a von Neumann algebra $\M$ of bounded operators, 
%
whereas the operators one wishes to integrate are often unbounded: generators of evolution equations, densities and moments in noncommutative probability, and Hamiltonians. 
For a positive self-adjoint operator $A$ affiliated with $\M$ there are several natural candidates for ``$\tau(A)$'': the spectral integral $\int\lambda\,d\tau(E_A(\lambda))$; the sharp truncations $\sup_R\tau(AE_A([0,R]))$; the clippings $\sup_R\tau(\min(A,R\id))$; the Yosida regularizations $\sup_h\tau(A(\id+hA)^{-1})$; the semigroup quotients $\lim_{t\downarrow0}t^{-1}\tau(\id-e^{-tA})$; and the supremum of $\tau(a)$ over $a\in\M_+$ with $a\le A$. If $\tau$ is normal, all of them coincide (Lemma~\ref{lem:normal-spectral}, Propositions~\ref{prop:clip} and~\ref{prop:exp}, and~\eqref{eq:inner}), and the classical theory of noncommutative integration applies~\cite{Nel74,FK86}. Our objective is to \emph{recover and approximate} the extended trace by evaluating the original bounded trace on explicitly constructed resolvents, semigroups and Chernoff products, rather than to start from a spectral integral valid only for normal traces. The question becomes substantive when $\tau$ is not normal.

\smallskip\noindent\textit{Quantum-mechanical interpretation.}  If $\tau$ is a tracial state and $H\ge0$ is affiliated with $\M$, then $\ty(H)$ is the mean energy in the tracial state, that is, at infinite temperature, and by Proposition~\ref{prop:exp}
\[
\ty(H)=\lim_{\beta\downarrow0}\frac{Z(0)-Z(\beta)}{\beta}=-Z'(0^+),
\]
for every trace (normal or not), where $Z(\beta)=\tau(e^{-\beta H})$ is the partition function of the system with Hamiltonian $H$, and $\beta$ the thermodynamic beta which is inversely proportional to the temperature. Hence, the mean energy equals minus the derivative of the partition function at infinite temperature. 

\smallskip\noindent\textit{Why the Yosida extension.}
We take $\ty(A)=\sup_h\tau(A(\id+hA)^{-1})$ as the definition (Definition~\ref{def:Y}). It coincides with the clipped extension (Proposition~\ref{prop:clip}) and the semigroup limit (Proposition~\ref{prop:exp}), and for every trace
\[
\ty(A)=\sup\{\tau(a):\ a\in\M_+,\ a\le A\}
\]
by~\eqref{eq:inner}. The last identity makes $\ty$ canonical: since $\ty$ is monotone, it is the \emph{smallest} extension of $\tau$ that is monotone for the order of quadratic forms. Minimality is the right normalization: assigning $+\infty$ to every unbounded operator is also monotone and additive, but is useless. The resolvent formulation is also the one available in evolution equations: $A(\id+hA)^{-1}=h^{-1}\bigl(\id-(\id+hA)^{-1}\bigr)$ is the Yosida approximation of $A$, and evaluating it requires one implicit Euler step, the solution of $(\id+hA)u=f$.

\smallskip\noindent\textit{Why additivity is a crux.}
For this inner extension superadditivity is formal: if $a\le A$ and $b\le B$, then $a+b\le A\splus B$, so $\ty(A)+\ty(B)\le\ty(A\splus B)$. Subadditivity would follow in the same way from a Riesz decomposition, that is, from writing every $c\in\M_+$ with $c\le A\splus B$ as $c=a+b$ with $0\le a\le A$ and $0\le b\le B$. Such a decomposition exists in abelian algebras, but not for noncommuting operators. In $M_2(\mathbb C)$ (with the standard basis $e_1=(1,0),e_2=(0,1)$ of $\mathbb C^2$, matrix units $e_{ii}=\operatorname{diag}(\delta_{1,i},\delta_{2,i})$), the projection $c$ onto $\mathbb C(e_1+e_2)$ satisfies $c\le e_{11}+e_{22}=\id$, while every $a$ with $0\le a\le e_{11}$ and every $b$ with $0\le b\le e_{22}$ are diagonal. The proof of Theorem~\ref{thm:additivity} replaces the missing decomposition by the trace property. It writes $f_h(A\splus B)=K_A^*K_A+K_B^*K_B$ with $K_AK_A^*\le f_h(A)$ and $K_BK_B^*\le f_h(B)$, and uses $\tau(K^*K)=\tau(KK^*)$. 


Additivity was already a difficult step for bounded operators. Murray and von Neumann defined the trace on a factor of type $\mathrm{II}_1$ through the dimension function~\cite{MvN36}; uniqueness was comparatively easy, but additivity was established only in~\cite{MvN37}. In~\cite{MvN37} they also showed that the closed densely defined operators affiliated with such a factor admit an addition and a multiplication. For every finite von Neumann algebra these operations make $\Aff(\M)$ a $*$-algebra; see also~\cite{Liu11,KL14}. 
Theorem~\ref{thmA} proves additivity for this algebra of unbounded operators and for arbitrary, possibly non-normal, traces.

\smallskip\noindent\textit{Approximation.}
In evolution equations the semigroup $e^{-tA}$ is rarely explicit, and one computes Chernoff products $F(t/n)^n$ of a simpler family with $F(0)=\id$ and $F'(0)=-A$~\cite{Che68}. Two questions arise. First, must the generator be affiliated with $\M$? Theorem~\ref{thmC}(i) shows that affiliation is forced: if $F(h)\in\M$ for all $h$, the semigroup lies in $\M$ and its generator is affiliated with $\M$. Second, can $\ty(A)$ be computed from $F$? For the resolvent family and for the semigroup itself it can, and more generally two-sided resolvent comparison suffices, with a limit independent of the approximation (Theorem~\ref{thm:sandwich}). Chernoff tangency alone does not suffice. The one-step trace defects $h^{-1}\tau(\id-F(h))$ may fail to converge to $\ty(A)$ even for a faithful normal tracial state (Proposition~\ref{prop:counter2}), because tangency is a pointwise statement on a core, whereas a trace integrates over the whole space. For normal traces the iterated limit $\lim_{t\downarrow0}\lim_{n\to\infty}t^{-1}\tau(\id-F(t/n)^n)$ does recover $\ty(A)$ (Corollary~\ref{cor:tracechernoff}).

\subsection{Setting}
Let $H$ be a complex separable Hilbert space and $\M\subseteq\mathcal B(H)$ a von Neumann algebra, where $\mathcal B(H)$ is the algebra of bounded linear{ operators on} $H$. 

The aim is to define a finite linear extension on a suitable class of unbounded operators, using bounded regularizations. Three questions arise: (i) whether the unbounded operators need to be affiliated; (ii) whether Chernoff approximations can be used to compute effectively traces; and (iii) whether the regularized positive functional is additive and extends to a linear trace. {For a positive self-adjoint operator $A$ affiliated with $\M$ we use the Yosida regularization (Definition~\ref{def:Y})
\begin{equation}\label{eq:Yintro}
\ty(A)=\lim_{h\downarrow0}\tau\bigl(A(\id+hA)^{-1}\bigr)=\sup_{h>0}\tau\bigl(A(\id+hA)^{-1}\bigr)\in[0,\infty].
\end{equation}
}

\subsection{Main results}
For finite $\M$, $\Aff(\M)$ denotes the $*$-algebra of operators affiliated with $\M$ (\S\ref{sec:formsum}), and $A\splus B$ denotes the form sum{; for finite $\M$ it is the sum in $\Aff(\M)$ (Lemma~\ref{lem:formsumequals})}.

\begin{mainthm}[Positive additivity and linear extension]\label{thmA}
Let $\tau$ be a trace on a finite von Neumann algebra $\M$. Let $A,B\ge0$ be self-adjoint operators affiliated with $\M$. 
 Then
\[
\ty(A\splus B)=\ty(A)+\ty(B)\qquad\text{in }[0,\infty].
\]

 Moreover, $\ty$ extends uniquely to a linear functional $\widetilde\tau$ on the $\mathbb{C}$-linear span $\mathscr D_\tau$ of $\{A\in\Aff(\M):A\ge0,\ \ty(A)<\infty\}$; $\widetilde\tau$ is positive and agrees with $\tau$ on $\M$.

If, in addition, the trace $\tau$ is normal then
the set function $\nu_A(B)=\tau(E_A(B))$ is well-defined, is a finite positive Borel measure on $[0,\infty)$, 
\[
\ty(A)=\int_{[0,\infty)}\lambda\,d\nu_A(\lambda)=\sup_{R>0}\tau(AE_A([0,R]))=\sup_{R>0}\tau(\min(A,R\id))\in [0,\infty].
\]
\end{mainthm}

\begin{mainthm}[Quadratic structure]\label{thmB}
Let $\M$ be finite and $\tau$ a trace on $\M$. Then $\mathscr L^2_Y=\{T\in\Aff(\M):\ty(T^*T)<\infty\}$ is a vector space. The form $\langle T,S\rangle_Y=\widetilde\tau(T^*S)$ is a positive semidefinite Hermitian form on it, and it is definite if and only if $\tau$ is faithful. If $\tau$ is faithful and normal, then $\mathscr L^2_Y=L^2(\M,\tau)$; 
in particular, $\mathscr L^2_Y$ is complete.
\end{mainthm}

\begin{mainthm}[Chernoff families]\label{thmC}
\begin{enumerate}[label=\textup{(\roman*)}]
\item Let $L$ generate a $C_0$-semigroup $T$ on $H$, and let $F(h)\in\M$ for $h\ge0$. If $F(t/n)^n\to T(t)$ strongly for every $t\ge0$, then $L$ is affiliated with $\M$.

\item There are positive contractive Chernoff families $F$ for the zero generator, in commutative von Neumann algebras and with $F(h)x-x=o(h)$ on a dense subspace, and a faithful normal tracial state $\tau$ such that $h^{-1}\tau(\id-F(h))$ does not converge to $\ty(0)=0$. 
\item Let $A\ge0$ be self-adjoint and affiliated with $\M$, let $\tau$ be a trace, and let $F(h)\in\M$. Suppose $A_{\alpha(h)}\le h^{-1}(\id-F(h))\le A_{\beta(h)}$, where $A_s=A(\id+sA)^{-1}$ and $\alpha(h),\beta(h)\to0^+$. Then $h^{-1}\tau(\id-F(h))\to\ty(A)$ in $[0,\infty]$. This applies to $F(h)=(\id+hA)^{-1}$ and to $F(h)=e^{-hA}$.
\item Let $A\geq0$ be self-adjoint affiliated with $\M$ and let $\tau$ be a normal trace. Suppose $F(h)\in\M$ and, for each fixed $t>0$, $F(t/n)^n\to e^{-tA}$ strongly as $n\to\infty$. Then
$$
\tau_Y(A)=\lim_{t\downarrow0}\lim_{n\to\infty}\frac1t\tau\bigl(\id-F(t/n)^n\bigr)\in[0,\infty].
$$
\end{enumerate}
\end{mainthm}

\begin{mainthm}[Integrability for arbitrary traces]\label{thmD}
Let $\M$ be finite, let $\tau$ be a trace on $\M$, and set $\mathscr L^1_Y=\{T\in\Aff(\M):\ty(|T|)<\infty\}$. For all $S,T\in\Aff(\M)$ and $X\in\M$,
\begin{gather*}
\ty(|S+T|)^{2/3}\le\ty(|S|)^{2/3}+\ty(|T|)^{2/3},\\
\ty(|XT|)\le\|X\|\,\ty(|T|),\qquad
\ty(|TX|)\le\|X\|\,\ty(|T|),
\end{gather*}
and
\begin{equation}\label{eq:AbsSum}
\tfrac12\bigl(\ty(|S+T|)+\ty(|S-T|)\bigr)\le\ty(|S|)+\ty(|T|).
\end{equation}
Consequently $\mathscr L^1_Y$ coincides with the domain $\mathscr D_\tau$ of $\widetilde\tau$, and it is a $\M$-bimodule. For $T\in\mathscr L^1_Y$ and $X\in\M$,
\begin{equation}\label{eq:tiltauatoY}
|\widetilde\tau(T)|\le\sqrt2\,\ty(|T|),\qquad \widetilde\tau(XT)=\widetilde\tau(TX).
\end{equation}
\end{mainthm}

%

\begin{remark} Then, $\mathscr L^1_Y$ has the following characterization: for any $T\in\Aff(\M)$,
 $$T\in\mathscr L^1_Y \Longleftrightarrow\ty((\operatorname{Re}T)_\pm),\ty((\operatorname{Im}T)_\pm)<\infty.$$
\end{remark}

\subsection{Ideas of proofs}

The proof of Theorem~\ref{thmA} takes no limit inside $\tau$. It is based on two comparisons between bounded operators at fixed $h>0$. The first is the operator inequality $f_{2h}(A)+f_{2h}(B)\le f_h(A\splus B)$, which follows from a variational formula for the inverse of a form operator. The second is a factorization $f_h(A\splus B)=K_A^*K_A+K_B^*K_B$ with $K_AK_A^*\le f_h(A)$ and $K_BK_B^*\le f_h(B)$; traciality turns it into the reverse inequality for traces.

Theorem~\ref{thmD} rests on two observations. First, $\ty(|T|)$ is the supremum of $\tau(c)$ over all $c\in\M_+$ with $\|cx\|\le\|Tx\|$ for $x\in\Dom T$ (Corollary~\ref{cor:var}); thus operator inequalities for $T^*T$ yield inequalities for $\ty(|T|)$. Second, integrating Theorem~\ref{thm:additivity}(b) over $h$ gives the square-root inequality $\ty((X+Y)^{1/2})\le\ty(X^{1/2})+\ty(Y^{1/2})$ (Proposition~\ref{prop:sqrt}). Combined with the elementary inequality $|S+T|^2\le(1+t)|S|^2+(1+t^{-1})|T|^2$, they imply $\ty(|S+T|)\le\sqrt{1+t}\,\ty(|S|)+\sqrt{1+t^{-1}}\,\ty(|T|)$ for every $t>0$. Optimizing over $t$ yields the exponent $2/3$, and in particular $\ty(|S+T|)\le\sqrt2\,\bigl(\ty(|S|)+\ty(|T|)\bigr)$. 

Theorem~\ref{thmA} combines Lemma~\ref{lem:formsumequals}, Lemma~\ref{lem:normal-spectral}, Theorem~\ref{thm:additivity} and Theorem~\ref{thm:linear}. Theorem~\ref{thmB} combines Theorems~\ref{thm:L2general} and~\ref{thm:L2normal}. Theorem~\ref{thmC} collects Theorem~\ref{thm:aff}, 
Proposition~\ref{prop:counter2}, Theorem~\ref{thm:sandwich}, Proposition~\ref{prop:exp} and Corollary~\ref{cor:tracechernoff}. {Theorem~\ref{thmD} combines Theorems~\ref{thm:L1quasi} and~\ref{thm:L1} and Corollary~\ref{cor:average}.}

\subsection{Organization and notation}
Section~\ref{sec:prelim} consists of some preliminaries and facts on affiliated operators and form sums. Section~\ref{sec:yosida} introduces $\ty$ and compares it with other truncations. Section~\ref{sec:chernoff} treats Chernoff families. Section~\ref{sec:additivity} proves positive additivity. Section~\ref{sec:linear} constructs the linear extension, the quadratic structure and contains the integrability results of Theorem~\ref{thmD}. Section~\ref{sec:normal} specializes to faithful normal traces. 

Throughout, $f_h(\lambda)=\lambda/(1+h\lambda)$ and $g_R(\lambda)=\min(\lambda,R)$ for $\lambda\ge0$ and $h,R>0$. We write $A_h=f_h(A)$, and $E_A$ for the spectral measure of a self-adjoint operator $A$. Order relations between unbounded positive operators are taken in the sense of quadratic forms. 

\section{Preliminaries}\label{sec:prelim}
\subsection{Von Neumann algebras and traces}\label{sec:vNa}
Throughout, $H$ is a separable complex Hilbert space and $\mathcal B(H)$ is the algebra of bounded operators on $H$. A { sequence $(X_n)$} in $\mathcal B(H)$ converges \emph{strongly} to $X\in\mathcal B(H)$ if { $X_nx\to Xx$} for every $x\in H$.

\begin{definition}[Von Neumann algebra]\label{def:vNa}
A \emph{von Neumann algebra} on $H$ is a subalgebra $\M\subseteq\mathcal B(H)$ such that
\begin{enumerate}[label=\textup{(\roman*)}]
\item $\id\in\M$;
\item $X\in\M$ implies $X^*\in\M$;
\item $\M$ is {sequentially} strongly closed: if $X_n\in\M$ and $X_nx\to Xx$ for every $x\in H$, then $X\in\M$.
\end{enumerate}
The \emph{commutant} of a set $\mathcal S\subseteq\mathcal B(H)$ is $\mathcal S'=\{B\in\mathcal B(H):BX=XB\text{ for all }X\in\mathcal S\}$. By von Neumann's bicommutant theorem, a unital $*$-subalgebra $\M$ of $\mathcal B(H)$ satisfies~(iii) if and only if $\M=\M''$~\cite[Ch.~5]{KR83}. We denote by $\M_+$ the set of positive elements of $\M$ and $\Z=\M\cap\M'$ for its centre.
\end{definition}

\begin{definition}[Positive functionals and traces]\label{def:trace}
Let $\M$ be a von Neumann algebra. A linear functional $\varphi:\M\to\mathbb C$ is \emph{positive} if $\varphi(X)\ge0$ for all $X\in\M_+$.\footnote{Then $|\varphi(X)|\le\varphi(\id)\|X\|$, so $\varphi$ is norm continuous with $\|\varphi\|=\varphi(\id)$.}  A \emph{trace} on $\M$ is a positive linear functional $\tau$ such that
\[
\tau(XY)=\tau(YX)\qquad(X,Y\in\M).
\]
A trace with $\tau(\id)=1$ is called a \emph{tracial state}.
\end{definition}

\begin{definition}[Faithful and normal functionals]\label{def:faithnormal}
Let $\varphi$ be a positive linear functional on a von Neumann algebra $\M$.
\begin{enumerate}[label=\textup{(\roman*)}]
\item $\varphi$ is \emph{faithful} if $X\in\M_+$ and $\varphi(X)=0$ imply $X=0$. Equivalently, $\varphi(X^*X)=0$ implies $X=0$.
\item $\varphi$ is \emph{normal (or $\sigma$-normal)} if $\varphi(X)=\sup_n\varphi(X_n)$ whenever $(X_n)$ is a bounded increasing sequence in $\M_+$ with strong limit $X$.
\end{enumerate}
\end{definition}


\subsection{Affiliated operators}
\begin{definition}[Affiliation]\label{def:commutant}
 A densely defined closed operator $L$ is \emph{affiliated} with $\M$, written $L\eta\M$, if for every unitary $U\in\M'$ one has $U\Dom(L)=\Dom(L)$ and $LUx=ULx$ for $x\in\Dom(L)$ {or, equivalently, $U^*LU=L$ for every unitary $U\in\M'$.}
\end{definition}

\begin{lemma}[One resolvent characterizes affiliation]\label{lem:one-resolvent}
Let $L$ be densely defined and closed, with $\rho(L)\neq\varnothing$. Then the following are equivalent:
\begin{enumerate}[label=\textup{(\roman*)}]
\item $L\eta\M$;
\item $(\lambda_0\id-L)^{-1}\in\M$ for some $\lambda_0\in\rho(L)$;
\item $(\lambda\id-L)^{-1}\in\M$ for every $\lambda\in\rho(L)$.
\end{enumerate}
\end{lemma}
\begin{proof}
(i)$\Rightarrow$(iii). Let $L\eta\M$. Then, for every unitary $U\in\M'$ and $\lambda\in\rho(L)$, $U^*(\lambda\id-L)^{-1}U=(U^*(\lambda\id-L)U)^{-1}=(\lambda\id-L)^{-1}$ i.e. $U$ commutes with the bounded resolvent. Since the unitaries linearly span $\M'$ then, by the bicommutant theorem, $(\lambda\id-L)^{-1}\in\M''=\M$.

 (iii)$\Rightarrow$(ii). Obvious.
 
 (ii)$\Rightarrow$(i). Let $R=(\lambda_0\id-L)^{-1}\in\M$. For unitary $U\in\M'$, $UR=RU$. Since $\Ran R=\Dom(L)$, both $U$ and $U^*$ preserve $\Dom(L)$. For $x=Ry\in\Dom(L)$,
\[
(\lambda_0\id-L)Ux=(\lambda_0\id-L)RUy=Uy=U(\lambda_0\id-L)x.
\]
Hence $LUx=ULx$; thus (ii)$\Rightarrow$(i).
\end{proof}

\begin{remark}[Where does the resolvent exist?]\label{rem:resolvent}
The resolvent is defined only for $\lambda\in\rho(L)$, not at arbitrary complex values. For $A\ge0$ self-adjoint, $(\lambda\id+A)^{-1}$ exists for each real $\lambda>0$, with norm at most $1/\lambda$. 
\end{remark}
\begin{lemma}[{\cite[Proposition~I.5.5 \& Theorem~II.1.10]{EN00}}]\label{lem:expostab}
If $L$ generates a $C_0$-semigroup then there exists $C\geq 1,\ \omega\in \mathbb{R}$ such that $\|e^{tL}\|\le Ce^{\omega t}$ for every $t\geq 0$; in particular, $\{\operatorname{Re}\lambda>\omega\}\subset\rho(L)$ and
\[
(\lambda\id-L)^{-1}x=\int_0^\infty e^{-\lambda t}e^{tL}x\,dt\qquad(\operatorname{Re}\lambda>\omega),
\]
a convergent vector-valued Bochner integral for every $x\in H$.
\end{lemma}

\begin{proposition}[Spectral projections and Borel functional calculus]\label{prop:spec-aff}
For a self-adjoint operator $A$, the following are equivalent:
\[
A\eta\M\quad\Longleftrightarrow\quad E_A(B)\in\M\quad\text{for every Borel }B\subseteq\mathbb R,
\]
where $E_A$ is the projection-valued measure map on $\mathbb R$ of $A$.

If these conditions hold, then $f(A)\in\M$ for every bounded Borel $f:\mathbb R\to\mathbb C${, and $f(A)\eta\M$ for every Borel $f:\mathbb R\to\mathbb C$. If $f$ is real-valued, then $f(A)$ is self-adjoint, and if $f\ge0$, then $f(A)\ge0$}.\footnote{Examples are $A^{1/2}$ for $A\ge0$, and $|A|$.}
\end{proposition}
{The proof is given in Appendix~\ref{app:specaff}.}

\subsection{The affiliated algebra of a finite von Neumann algebra and form sums}\label{sec:formsum}
{Recall that $\M$ is \emph{finite} if every isometry in $\M$ is unitary.}
The affiliated operators of a finite von Neumann algebra form a $*$-algebra, denoted $\Aff(\M)$, with addition and multiplication understood as closures of the corresponding densely defined operations~\cite{Liu11,MvN36,MvN37}. 
Following~\cite{Liu11}, we write $S\splus T$ and $S\sdot T$ for these closures. We shall simply write $S+T$ and $ST$ when there is no confusion. In a general von Neumann algebra, arbitrary sums of positive affiliated operators demand more careful quadratic-form terminology. {For self-adjoint $A,B\ge0$ with $\Dom(A^{1/2})\cap\Dom(B^{1/2})$ dense, the \emph{form sum} $A\splus B$ is the positive self-adjoint operator associated with the closed form $t[y]=\|A^{1/2}y\|^2+\|B^{1/2}y\|^2$ on $\Dom(A^{1/2})\cap\Dom(B^{1/2})$; the corresponding{ sesquilinear} form is denoted by $t(y,x)$.} 

\begin{lemma}[Affiliated sum equals form sum]\label{lem:formsumequals}
Let $\M$ be finite and let $A,B\ge0$ be self-adjoint operators affiliated with $\M$. Then $\mathcal D:=\Dom(A^{1/2})\cap\Dom(B^{1/2})$ is dense, the form sum $A\splus B${, whose form domain is $\mathcal D$,} is affiliated with $\M${, and it equals the closure of $A+B$ on $\Dom A\cap\Dom B$, that is, the sum in the $*$-algebra $\Aff(\M)$}.
\end{lemma}
{The proof is given in Appendix~\ref{app:formsum}.}

\begin{lemma}[Positive elements of $\Aff(\M)$]\label{lem:posaff}

Let $\M$ be finite, and let $Z\in\Aff(\M)$ satisfy $\langle Zx,x\rangle\ge0$ for all $x\in\Dom Z$. Then $Z$ is self-adjoint and $\sigma(Z)\subseteq[0,\infty)$. We write $\Aff(\M)_+$ for the set of such $Z$, and, for self-adjoint $W,Z\in\Aff(\M)$, $W\le Z$ if $Z-W\in\Aff(\M)_+$.\footnote{Observe that for any von Neumann algebra (not necessarily finite), for any  $A\in\M_+$ we have $X=A^{1/2}\in\M_+$ (by Proposition~\ref{prop:spec-aff}) and $A=X^*X$.
} Moreover:
\begin{enumerate}[label=\textup{(\roman*)},leftmargin=*]
\item $T^*T\in\Aff(\M)_+$ for every $T\in\Aff(\M)$, and $T^*T$ is the product $T^*\sdot T$ in $\Aff(\M)$;
\item $\Aff(\M)_+$ is closed under sums and under multiplication by nonnegative scalars, and $\le$ is transitive;
\item if $W,Z\in\Aff(\M)_+$ and $W\le Z$, then $W\le Z$ in the sense of quadratic forms; if $W\in\M_+$, the converse holds as well.
\end{enumerate}
\end{lemma}
The proof is given in Appendix~\ref{app:posaff}.

\section{The Yosida extension}\label{sec:yosida}
For $A\ge0$, $A\eta\M$, define
\[
A_h:=A(\id+hA)^{-1}=f_h(A),\qquad f_h(\lambda)=\frac{\lambda}{1+h\lambda},\qquad h>0.
\]

\begin{lemma}[Bounded monotone Yosida approximants]\label{lem:yosida-monotone}
For $A\ge0$ self-adjoint and $h>0$, $A_h\in\M$ and $0\le A_h\le h^{-1}\id$. For $0<k<h$, $A_h\le A_k$. Furthermore, $A_hx\to Ax$ in $H$ for every $x\in\Dom(A)$, and $\langle A_hx,x\rangle\uparrow\|A^{1/2}x\|^2$ for every $x\in\Dom(A^{1/2})$.
\end{lemma}
\begin{proof}
The scalar function $f_h(\lambda)=\lambda/(1+h\lambda)$ is smooth on $[0,\infty)$ and satisfies $0\le f_h(\lambda)\le 1/h$ and $f_h(\lambda)\le f_k(\lambda)$ if $k<h$. Then, Proposition~\ref{prop:spec-aff} yields $A_h\in\M$, and the spectral functional calculus gives the boundedness and order claims. If $x\in\Dom(A)$, with $\mu_x(B)=\|E_A(B)x\|^2$, then
\[
\|(A_h-A)x\|^2=\int_{[0,\infty)}\left(\frac{h\lambda^2}{1+h\lambda}\right)^2d\mu_x(\lambda)\to0,
\]
by dominated convergence: the integrand tends pointwise to zero and is dominated by $\lambda^2$, integrable precisely for $x\in\Dom(A)$. If $x\in\Dom(A^{1/2})$, scalar monotone convergence gives
\[
\langle A_hx,x\rangle=\int f_h(\lambda)\,d\mu_x(\lambda)\uparrow\int\lambda\,d\mu_x(\lambda)=\|A^{1/2}x\|^2.
\]
\end{proof}

\begin{remark}\label{rmk:MonCon}
Introduce the extended positive quadratic form of $A\geq 0$:
\[
q_A[x]=\begin{cases}\|A^{1/2}x\|^2,&x\in\Dom(A^{1/2}),\\ +\infty,&x\notin\Dom(A^{1/2}),\end{cases}
\]
and set $q_h[x]=q_{A_h}[x]$.  
Then $q_h[x]\leq q_k[x]$ whenever $0<k<h$, and $\sup_{h>0}q_h[x]=q_A[x]$ for every $x\in H$. Thus, Lemma~\ref{lem:yosida-monotone} implies monotone (increasing) convergence of \emph{extended quadratic forms} as $h\downarrow0$, which we {abbreviate as} $A_h\uparrow A$.
\end{remark}
\begin{definition}\label{def:Y}
{For a trace $\tau$ on $\M$ and $A\ge0$ self-adjoint and affiliated with $\M$, we extend the trace as:}\footnote{This is a prescribed \emph{regularization extension}. It makes sense even if $\tau$ is not normal.}
\begin{equation}\label{eq:Y}
\ty(A):=\sup_{h>0}\tau(A_h)=\lim_{h\downarrow0}\tau(A_h)\in[0,+\infty].
\end{equation}
\end{definition}

\begin{proposition}[Elementary properties]\label{prop:elementary}
For every {trace $\tau$}, the map $\ty$ is well-defined with values in $[0,\infty]$, agrees with $\tau$ on $\M_+$, is positively homogeneous, is monotone for the positive quadratic-form order, and is invariant under conjugation by unitaries in $\M$.
\end{proposition}
\begin{proof}
\emph{Existence:} By Lemma~\ref{lem:yosida-monotone}, $0\le A_h\le A_k$ if $h>k>0$, so positivity gives $0\le\tau(A_h)\le\tau(A_k)$. Therefore $\ty(A)=\sup_h\tau(A_h)$ exists as an extended positive number.

\emph{Agreement:} If $A\in\M_+$, then $A-A_h=hA^2(\id+hA)^{-1}$ and $\|A-A_h\|\le h\|A\|^2$, since $\|(\id+hA)^{-1}\|\le1$. Norm continuity of $\tau$ yields $\ty(A)=\tau(A)$.

\emph{Homogeneity:} For $c>0$, $(cA)_h=cA_{ch}$, hence $\ty(cA)=c\sup_{s>0}\tau(A_s)=c\ty(A)$; the case $c=0$ is immediate (with the standard convention $0\cdot\infty=0$).

\emph{Monotonicity:} If $0\le A\le B$ in closed quadratic-form order, the inverse order relation gives $(\id+hB)^{-1}\le(\id+hA)^{-1}$. Subtracting from $\id$ and dividing by $h$ gives $A_h\le B_h$, so $\ty(A)\le\ty(B)$.

\emph{Unitary invariance:} If $U\in\M$ is unitary, Borel functional calculus yields $(UAU^*)_h=UA_hU^*$. Cyclicity for bounded elements gives $\tau((UAU^*)_h)=\tau(A_h)$ for every $h$ and taking the suprema yields the claim.
\end{proof}

\begin{lemma}[Comparison of Yosida and clipped regularizations]\label{lem:clip-comparison}
For $R,h>0$ and $\lambda\geq0$,
\[
f_{1/R}(\lambda)\leq g_R(\lambda),\qquad f_h(\lambda)\leq g_{1/h}(\lambda),\qquad g_R(\lambda)\leq(1+hR)f_h(\lambda).
\]
Consequently the corresponding inequalities hold for the bounded operators $f_h(A)$ and $g_R(A)$.
\end{lemma}
\begin{proof}
For $\lambda\leq R$, $g_R(\lambda)=\lambda$ and $f_{1/R}(\lambda)=\lambda/(1+\lambda/R)\leq\lambda$. For $\lambda\geq R$, $g_R(\lambda)=R$ and $f_{1/R}(\lambda)=R\lambda/(R+\lambda)\leq R$. Replacing $R$ with $1/h$ proves $f_h\leq g_{1/h}$. For the last inequality, if $0<\lambda\leq R$ then $g_R(\lambda)/f_h(\lambda)=1+h\lambda\leq1+hR$ (the case $\lambda=0$ is immediate). If $\lambda\geq R$,
\[
(1+hR)f_h(\lambda)-R=\frac{(1+hR)\lambda-R(1+h\lambda)}{1+h\lambda}=\frac{\lambda-R}{1+h\lambda}\geq0.
\]
All comparisons follow pointwise; bounded Borel functional calculus preserves the order of real functions, yielding the operator versions.
\end{proof}

\begin{proposition}[Clipped spectral representation]\label{prop:clip} For every trace $\tau$ and every self-adjoint $A\ge0$ affiliated with $\M$,\footnote{Without faithfulness and normality.}
\begin{equation}\label{eq:clip}
\ty(A)=\sup_{R>0}\tau(\min(A,R\id)).
\end{equation}
\end{proposition}
\begin{proof}
By Lemma~\ref{lem:clip-comparison}, $f_h\le g_{1/h}$; thus $\sup_h\tau(A_h)\le\sup_R\tau(g_R(A))$. The {third inequality} gives $\tau(g_R(A))\le(1+hR)\tau(A_h)$; letting $h\downarrow0$ proves the reverse inequality for each $R$.
\end{proof}

\begin{lemma}[Spectral integral of a normal positive trace]\label{lem:normal-spectral}
Let $\tau$ be a positive normal tracial functional on $\M$, and let $A\ge0$ be self-adjoint and affiliated with $\M$. The set function $\nu_A(B)=\tau(E_A(B))$ is a finite positive Borel measure on $[0,\infty)$, and, in $[0,\infty]$,
\[
\ty(A)=\int_{[0,\infty)}\lambda\,d\nu_A(\lambda)=\sup_{R>0}\tau(AE_A([0,R]))=\sup_{R>0}\tau(\min(A,R\id)).
\]
\end{lemma}
\begin{proof}
The spectral projections lie in $\M$ by Proposition~\ref{prop:spec-aff}. If $B_j$ are pairwise disjoint Borel sets, then $P_N=\sum_{j=1}^NE_A(B_j)\uparrow E_A(\cup_{j\ge1}B_j)$ strongly. Normality yields
\[
\nu_A\Bigl(\bigcup_{j\ge1}B_j\Bigr)=\lim_N\tau(P_N)=\sum_{j\ge1}\nu_A(B_j),
\]
so $\nu_A$ is countably additive; $\nu_A([0,\infty))=\tau(\id)<\infty$. For a bounded nonnegative simple Borel function $f=\sum_j a_j\mathbf1_{B_j}$, functional calculus and linearity give $\tau(f(A))=\sum_ja_j\tau(E_A(B_j))=\int f\,d\nu_A$. For a bounded nonnegative Borel $f$, approximate it by increasing simple functions $f_n\uparrow f$. Then $f_n(A)\uparrow f(A)$ strongly, so normality on positive increasing bounded sequences, together with scalar monotone convergence, gives
\[
\tau(f(A))=\lim_n\tau(f_n(A))=\lim_n\int f_n\,d\nu_A=\int f\,d\nu_A.
\]
Take $f=f_h$ to obtain $\tau(A_h)=\int f_h\,d\nu_A$. As $h\downarrow0$, $f_h(\lambda)\uparrow\lambda$, whence $\ty(A)=\int\lambda\,d\nu_A$ by scalar monotone convergence. Finally, the bounded functions $\lambda\mathbf1_{[0,R]}(\lambda)$ and $\min(\lambda,R)$ both increase pointwise to $\lambda$ as $R\to\infty$. Then, applying the same bounded Borel trace identity and monotone convergence, the two truncation formulas follow.
\end{proof}

\begin{remark}\label{rem:normalused}
Normality is used for countable additivity of $\nu_A$ and for passing through increasing spectral approximations. For non-normal traces the clipped formula (see also Proposition~\ref{prop:clip}) remains valid but  spectral cutoffs need not yield the same extension.
\end{remark}

\section{Chernoff families and trace defects}\label{sec:chernoff}

\subsection{\texorpdfstring{Chernoff products in $\M$ force affiliation}{Chernoff products in M force affiliation}}
\begin{theorem}[Automatic affiliation]\label{thm:aff}
Let $L$ generate a $C_0$-semigroup $T(t)=e^{tL}$ with $\|T(t)\|\le Ce^{\omega t}$ for $t\geq 0$, and suppose $F(h)\in\M$ has Chernoff products $F(t/n)^n\to T(t)$ strongly for every $t\ge0$ as $n\to\infty$. Then $L\eta\M$. In particular, the assertion applies to $L=-A$ for $A\ge0$ self-adjoint.
\end{theorem}
\begin{proof}
For each $t$ and $n$, the product $F(t/n)^n$ lies in $\M$. Therefore  $T(t)\in\M$ for every $t\ge0$ (see Definition~\ref{def:vNa}(iii)). Choose  $\lambda_0>\omega$. The Hille--Yosida resolvent formula gives $R=(\lambda_0\id-L)^{-1}$ and, for each $x\in H$,
\[
Rx=\int_0^\infty e^{-\lambda_0t}T(t)x\,dt,\qquad \int_0^\infty e^{-\lambda_0t}\|T(t)x\|\,dt\le\frac{C}{\lambda_0-\omega}\|x\|.
\]
For a unitary $U\in\M'$, $UT(t)=T(t)U$ for every $t$. The map $t\mapsto e^{-\lambda_0t}T(t)x$ is strongly continuous and norm-integrable by the preceding estimate. For a bounded operator $U$ and $x\in H$, letting $v(t)=e^{-\lambda_0t}T(t)x$, it follows
\begin{align*}
\Bigl\langle U\!\int_0^\infty v(t)dt,\;y\Bigr\rangle=\Bigl\langle \int_0^\infty v(t)dt,\;U^*y\Bigr\rangle&=\int_0^\infty\langle v(t),U^*y\rangle\,dt\\
 &=\int_0^\infty\langle Uv(t),y\rangle\,dt=\Bigl\langle\int_0^\infty Uv(t)dt,\;y\Bigr\rangle ,
\end{align*}
for every $y\in H$, i.e.
\[
URx=\int_0^\infty e^{-\lambda_0t}UT(t)x\,dt=\int_0^\infty e^{-\lambda_0t}T(t)Ux\,dt=RUx.
\]
Thus $UR=RU$, and $R\in\M''=\M$. Lemma~\ref{lem:one-resolvent} now implies $L\eta\M$.
\end{proof}

\begin{corollary}[Generators of semigroups in $\M$]\label{cor:gen}
{Let $L$ generate a $C_0$-semigroup $T_L$ with $T_L(t)\in\M$ for all $t\ge0$. Then $L\eta\M$. Every $B\in\M$ generates such a semigroup.}
\end{corollary}
\begin{proof}
{Apply Lemma~\ref{lem:expostab} and Theorem~\ref{thm:aff} with $F(h):=T_L(h)$, for which $F(t/n)^n=T_L(t)$ for every $n$. For $B\in\M$, the series $e^{tB}=\sum_{k\ge0}t^kB^k/k!$ converges in norm in $\M$, the map $t\mapsto e^{tB}$ is norm continuous, and $\|e^{tB}\|\le e^{t\|B\|}$.}
\end{proof}

\subsection{Chernoff tangency does not control trace defects}
The resolvent family $F_A(h)=(\id+hA)^{-1}$ is Chernoff equivalent to the contraction semigroup $e^{-tA}$:
\[
F_A(t/n)^n\xrightarrow[n\to\infty]{strongly}e^{-tA}.
\]
{For self-adjoint $A\ge0$ this follows from the spectral theorem and dominated convergence, since $\lim_n (1+t\lambda/n)^{-n}= e^{-t\lambda}$ boundedly on $[0,\infty)$, which is the Euler formula, a particular case of Chernoff's product formula~\cite{Che68},~\cite[Ch.~III]{EN00}.} Moreover $(\id-F_A(h))/h=A_h$, so its infinitesimal defect recovers $\ty$ by definition. However, arbitrary Chernoff tangency is too weak for trace $\ty$ recovery.

\begin{definition}\label{def:chernoff}
{Let $(L,D(L))$ be a closed operator which generates a $C_0$-semigroup $T$ on $H$. A \emph{Chernoff function for $T$ in $\M$} is a strongly continuous map $F:[0,\infty)\to\mathcal B(H)$ with $F(0)=\id$, $F'(0)=L$ on a core $\mathcal{D}$ of $L$ with closure $(L,D(L))$  and $\lim_n F(t/n)^n=T(t)$ strongly{,} locally {uniformly} in $t\ge0$. It is \emph{positive contractive} if $0\le F(h)\le\id$ for all $h\geq 0$.}
\end{definition}

\begin{proposition}[Counterexample with a normal trace]\label{prop:counter2}
{Let $\M=\ell^\infty$ on $\ell^2(\mathbb N)$ with the faithful normal tracial state $\tau(x)=\sum_{k\ge1}2^{-k}x_k$. There is a positive contractive Chernoff function $F$ for the zero generator in $\M$, with $F(h)x=x$ for all sufficiently small $h$ whenever $x$ is finitely supported, such that $2^m\tau(\id-F(2^{-m}))\ge\frac12$ for every $m\ge1$.}
\end{proposition}
\begin{proof}
 Choose $\psi\in C_c^\infty((1,4))$, $0\le\psi\le1$, $\psi(2)=1$, and let
\[
F(h)=\operatorname{diag}(1-\psi(2^kh))_{k\ge1},\qquad F(0)=\id.
\]
 {Let $x$ be finitely supported, say $x_k=0$ for $k>K$. Since $\operatorname{supp}\psi\subset(1,4)$, we have $\psi(2^kh)=0$ whenever $2^kh\le1$. Hence $F(h)x=x$ for $0\le h\le2^{-K}$. In particular, $F$ is tangent to the zero generator on this dense subspace, and $F(t/n)^nx=x$ as soon as $n\ge2^Kt$. To pass to arbitrary $x\in\ell^2$ we use the uniform bound $\|F(t/n)^n\|\le1$, which holds because $0\le F(h)\le\id$. Given $x\in\ell^2$ and $\varepsilon>0$, choose a finitely supported $y$, with $y_k=0$ for $k>K$, such that $\|x-y\|<\varepsilon$. Then, for $n\ge2^Kt$,
\[
\|F(t/n)^nx-x\|\le\|F(t/n)^n(x-y)\|+\|F(t/n)^ny-y\|+\|y-x\|<2\varepsilon .
\]
Hence $F(t/n)^n\to\id$ strongly. 
 The density step cannot be skipped..} For $h_m=2^{-m}$, the $k=m+1$ summand gives
\[
\frac1{h_m}\tau(\id-F(h_m))\ge\frac{2^{-(m+1)}}{2^{-m}}\psi(2)=\frac12.
\]
Therefore the trace defect does not converge to $\ty(0)=0$.
\end{proof}

\subsection{Two-sided resolvent comparison}
\begin{theorem}[Trace-defect sandwich]\label{thm:sandwich}
Let $A\ge0$ be affiliated, {let $\tau$ be a trace,} and let $F_A(h)\in\M$ satisfy $D_h=(\id-F_A(h))/h\in\M_+$. If $\alpha(h),\beta(h)>0$ tend to zero and
\[
A_{\alpha(h)}\le D_h\le A_{\beta(h)},
\]
then $\lim_{h\downarrow0}\tau(D_h)=\ty(A)$ in $[0,+\infty]$.
\end{theorem}
\begin{proof}
Apply positivity of $\tau$ to the two inequalities and pass to the extended monotone limits.
\end{proof}

\begin{proposition}[Exponential regularization]\label{prop:exp}
For every positive affiliated $A$ and every positive {trace $\tau$},
\[
\ty(A)=\lim_{h\downarrow0}\frac1h\tau(\id-e^{-hA}).
\]
\end{proposition}
\begin{proof}
{For $x\ge0$},
\[
\frac{x}{1+x}\le1-e^{-x}\le\frac{x}{1+x/2}.
\]
Functional calculus yields $A_h\le(\id-e^{-hA})/h\le A_{h/2}$. Then, apply Theorem~\ref{thm:sandwich}.
\end{proof}

\begin{problem}[Approximation-independence]\label{prob:indep}
Identify useful classes of positive contraction Chernoff families $F_A(h)\in\M$ for which
\[
\lim_{h\downarrow0}\tau\Bigl(\frac{\id-F_A(h)}h\Bigr)=\ty(A)
\]
for every positive affiliated $A$ and {every positive trace} $\tau$.
\end{problem}

\section{Positive additivity}\label{sec:additivity}
{
\subsection{Strongly commuting positive affiliated operators}
\begin{theorem}[Commuting additivity]\label{thm:commuting}
If $A,B\ge0$ are strongly commuting affiliated operators, then
\[
\ty(A\splus B)=\ty(A)+\ty(B)
\]
with equality in $[0,\infty]$.\footnote{It holds without requiring normality and faithfulness. Recall that $A\splus B$ denotes the closure of the sum on $\Dom A\cap\Dom B$, which is the self-adjoint operator $\int(a+b)\,dE(a,b)$ given by the joint spectral measure $E$ of $(A,B)$; it is affiliated with $\M$  by Lemma~\ref{lem:affcrit}. 
}
\end{theorem}
\begin{proof}
For $a,b\ge0$, one has
\[
g_{R/2}(a)+g_{R/2}(b)\le g_R(a+b)\le g_R(a)+g_R(b).
\]
Then, use the joint spectral calculus, apply $\tau$, and let $R\to\infty$ using Proposition~\ref{prop:clip}.
\end{proof}
}

\subsection{Positive additivity via form--resolvent comparison}\label{sec:formsumadd}
\begin{theorem}[Positive additivity without normality]\label{thm:additivity}
Let $\M$ be any von Neumann algebra and $\tau$ any positive linear functional on $\M$ with $\tau(XY)=\tau(YX)$. Let $A,B\ge0$ be self-adjoint and affiliated with $\M$, with $\mathcal D=\Dom(A^{1/2})\cap\Dom(B^{1/2})$ dense, and let $S:=A\splus B$. For every $h>0$:
\begin{enumerate}[label=\textup{(\alph*)}]
\item $f_{2h}(A)+f_{2h}(B)\le f_h(S)$ in $\M$;
\item $\tau(f_h(S))\le\tau(f_h(A))+\tau(f_h(B))$;
\item $\ty(A\splus B)=\ty(A)+\ty(B)$ in $[0,\infty]$.
\end{enumerate}
\end{theorem}
\begin{proof}
$S$ is affiliated with $\M$ by  Lemma~\ref{lem:formsumequals}, which does not use finiteness.

\emph{(a)} Let $T\ge\id$ be self-adjoint. Since $\sigma(T)\subseteq[1,\infty)$, the operator $T^{-1/2}$ is bounded and $\Dom(T^{-1/2})=H$. For $y\in\Dom(T^{1/2})$ one has $\langle x,y\rangle=\langle T^{-1/2}x,T^{1/2}y\rangle$, hence
\[
2\operatorname{Re}\langle x,y\rangle-\|T^{1/2}y\|^2=\|T^{-1/2}x\|^2-\|T^{1/2}y-T^{-1/2}x\|^2 .
\]
Taking $y=T^{-1}x$ gives the variational formula
\begin{equation}\label{eq:var}
\langle T^{-1}x,x\rangle=\sup_{y\in\Dom(T^{1/2})}\bigl(2\operatorname{Re}\langle x,y\rangle-\|T^{1/2}y\|^2\bigr).
\end{equation}
Apply it {with} $P=\id+2hA$, $Q=\id+2hB$ and $T=\id+hS$. For $y\in\mathcal D=\Dom(T^{1/2})$,
\[
\|T^{1/2}y\|^2=\|y\|^2+h\,t[y]=\tfrac12\|P^{1/2}y\|^2+\tfrac12\|Q^{1/2}y\|^2 .
\]
Hence
\[
\begin{aligned}
\langle T^{-1}x,x\rangle&=\sup_{y\in\mathcal D}\Bigl[\tfrac12\bigl(2\operatorname{Re}\langle x,y\rangle-\|P^{1/2}y\|^2\bigr)+\tfrac12\bigl(2\operatorname{Re}\langle x,y\rangle-\|Q^{1/2}y\|^2\bigr)\Bigr]\\
&\le\tfrac12\langle P^{-1}x,x\rangle+\tfrac12\langle Q^{-1}x,x\rangle .
\end{aligned}
\]
Since $f_h(S)=h^{-1}(\id-T^{-1})$ and $f_{2h}(A)=(2h)^{-1}(\id-P^{-1})$, (a) follows.

\emph{(b)} Let $G=(\id+hS)^{-1/2}\in\M$. Then $\Ran G=\Dom(S^{1/2})=\mathcal D\subseteq\Dom(A^{1/2})$, so $K_A:=A^{1/2}G$ is everywhere defined. It is bounded:
\[
\|K_Ax\|^2\le \langle S(\id+hS)^{-1}x,x\rangle\le h^{-1}\|x\|^2 .
\]
Moreover $K_A\in\M$ by Lemma~\ref{lem:affcrit}(ii).
 The same holds for $K_B:=B^{1/2}G$ likewise. For every $x\in \mathcal D$,
\[
\|K_Ax\|^2+\|K_Bx\|^2=\|S^{1/2}Gx\|^2=\langle f_h(S)x,x\rangle,
\]
hence $K_A^*K_A+K_B^*K_B=f_h(S)$.\footnote{For $x\in\Dom(A^{1/2})$ one checks $K_A^*x=GA^{1/2}x$.}  Since $A\le S$ then $(\id+hS)^{-1}\le(\id+hA)^{-1}$. Therefore
\[
\|K_A^*x\|^2=\langle(\id+hS)^{-1}A^{1/2}x,A^{1/2}x\rangle\le\langle(\id+hA)^{-1}A^{1/2}x,A^{1/2}x\rangle=\langle f_h(A)x,x\rangle .
\]
By density and continuity, $K_AK_A^*\le f_h(A)$, and similarly $K_BK_B^*\le f_h(B)$. Traciality and positivity give
\[
\tau(f_h(S))=\tau(K_AK_A^*)+\tau(K_BK_B^*)\le\tau(f_h(A))+\tau(f_h(B)).
\]

\emph{(c)} By (a), (b) and positivity,
\[
\tau(f_{2h}(A))+\tau(f_{2h}(B))\le\tau(f_h(S))\le\tau(f_h(A))+\tau(f_h(B)).
\]
As $h\downarrow0$ each $\tau(f_h(\cdot))$ increases to $\ty(\cdot)$ by {Definition~\ref{def:Y} and Lemma~\ref{lem:yosida-monotone}}, so both bounds tend to $\ty(A)+\ty(B)$.
\end{proof}

\begin{remark}\label{rem:additivity}
Theorem~\ref{thm:additivity} uses neither normality  nor faithfulness of $\tau$, and no limit is taken inside $\tau$.
\end{remark}

\begin{remark}[Alternative proof of superadditivity]
\label{rem:BH}
For every trace $\tau$ and every positive self-adjoint $A$ affiliated with $\M$,
\begin{equation}\label{eq:inner}
\ty(A)=\sup\{\tau(a):a\in\M_+,\ a\le A\},
\end{equation}
where $a\le A$ is meant in the sense of quadratic forms. Indeed, if $a\le A$, then $(\id+ha)^{-1}\ge(\id+hA)^{-1}$, so that $f_h(a)\le f_h(A)$. Since $\|a-f_h(a)\|\le h\|a\|^2$, we get $\tau(a)=\lim_{h\downarrow0}\tau(f_h(a))$. 
Therefore, 
$
\sup\{\tau(a):a\in\M_+,\ a\le A\}=\sup_{h>0}\tau(A_h)=\ty(A).
$
  If $a\le A$ and $b\le B$, then $a+b\le A\splus B$, so \eqref{eq:inner} gives $\ty(A)+\ty(B)\le\ty(A\splus B)$. 
\end{remark}

\begin{remark}[Relation to Bourin--Hiai]\label{rem:BHcredit}
Let $\M$ be finite. Since $H$ is separable, $\M$ carries a faithful normal tracial state $\tau_0$~\cite[Ch.~8]{KR86}, and every positive self-adjoint operator affiliated with $\M$ is $\tau_0$-measurable~\cite{FK86}. For such $A,B$ and every $\varepsilon>0$, Theorem~3.2 of Bourin and Hiai~\cite{BH15}, applied with the concave function $f_h$, provides unitaries $U,V\in\M$ such that
\[
f_h(A\splus B)\le Uf_h(A)U^*+Vf_h(B)V^*+\varepsilon\id .
\]
Applying an arbitrary trace $\tau$, normal or not, and letting $\varepsilon\to0$ yields Theorem~\ref{thm:additivity}(b) for finite $\M$. Our proof is direct 
 and applies to every von Neumann algebra.
\end{remark}

\begin{corollary}[A variational formula for $\ty(|T|)$]\label{cor:var}
Let $\tau$ be a trace on $\M$, and let $T$ be a closed, densely defined operator affiliated with $\M$. Then
\[
\ty(|T|)=\sup\bigl\{\tau(c):\ c\in\M_+,\ \|cx\|\le\|Tx\|\ \text{for all }x\in\Dom T\bigr\}.
\]
\end{corollary}
Observe that $|T|=(T^*T)^{1/2}$ is positive, self-adjoint and affiliated with $\M$ (Proposition~\ref{prop:spec-aff}), with $\Dom|T|=\Dom T$ and $\||T|x\|=\|Tx\|$. Hence the condition on $c$ reads $c^2\le T^*T$ i.e.
\[
\ty(|T|)=\sup\bigl\{\tau(c):\ c\in\M_+,\ c^2\le T^*T\bigr\}.
\]
\begin{proof}
If $c\in\M_+$ and $c^2\le|T|^2$, then $c\le|T|$ by Lemma~\ref{lem:heinz}, so $\tau(c)\le\ty(|T|)$ by~\eqref{eq:inner}. Conversely, $c=f_h(|T|)$ satisfies $c^2\le|T|^2$ because $0\le f_h(\lambda)\le\lambda$, and $\tau(f_h(|T|))\uparrow\ty(|T|)$ as $h\downarrow0$.
\end{proof}

\subsection{A square-root inequality}\label{sec:sqrt}
For $\lambda\ge0$,
\[
\pi\lambda^{1/2}=\int_0^\infty s^{-1/2}\,\frac{\lambda}{\lambda+s}\,ds=\int_0^\infty s^{-3/2}f_{1/s}(\lambda)\,ds ,
\]
so the square root is a superposition of the functions $f_h$. Integrating Theorem~\ref{thm:additivity}(b) over $h$ therefore gives an inequality of Rotfel'd type.

\begin{proposition}[Square-root subadditivity]\label{prop:sqrt}
Let $\tau$ be a trace on a von Neumann algebra $\M$, and let $X,Y\ge0$ be self-adjoint operators affiliated with $\M$ such that $\Dom(X^{1/2})\cap\Dom(Y^{1/2})$ is dense. Then
\[
\ty\bigl((X\splus Y)^{1/2}\bigr)\le\ty(X^{1/2})+\ty(Y^{1/2}).
\]
\end{proposition}
{
\begin{proof}[Idea of the proof]
For $0<\varepsilon<M$ let $g_{\varepsilon,M}(\lambda)=\pi^{-1}\int_\varepsilon^Ms^{-3/2}f_{1/s}(\lambda)\,ds$. For each positive self-adjoint $Z$ affiliated with $\M$, the operator $g_{\varepsilon,M}(Z)$ is a norm-convergent Riemann integral of the operators $s^{-3/2}f_{1/s}(Z)\in\M$, so $\tau$ may be applied under the integral sign. Theorem~\ref{thm:additivity}(b) with $h=1/s$ then gives $\tau(g_{\varepsilon,M}(X\splus Y))\le\tau(g_{\varepsilon,M}(X))+\tau(g_{\varepsilon,M}(Y))$. Since $g_{\varepsilon,M}(\lambda)\le\lambda^{1/2}$, the right-hand side is at most $\ty(X^{1/2})+\ty(Y^{1/2})$ by~\eqref{eq:inner}. The truncation is removed with scalar estimates only, namely $\min(\lambda^{1/2},R)\le g_{\varepsilon,M}(\lambda)+\pi^{-1}\bigl(2\varepsilon^{1/2}+2R^2M^{-1/2}\bigr)$, followed by Proposition~\ref{prop:clip}. The details are in Appendix~\ref{app:sqrt}.
\end{proof}
}

\begin{remark}\label{rem:powers}
The same argument, based on $\lambda^p=\frac{\sin p\pi}{\pi}\int_0^\infty s^{p-1}\frac{\lambda}{\lambda+s}\,ds$, gives $\ty((X\splus Y)^p)\le\ty(X^p)+\ty(Y^p)$ for $0<p<1$. For $p=1$ this is Theorem~\ref{thm:additivity}(c). 
Observe that these inequalities hold for traces, including the ones which are not normal: as in Theorem~\ref{thm:additivity}, no limit is taken inside $\tau$.
\end{remark}

\section{The linear extension and the quadratic structure}\label{sec:linear}
{Throughout this section $\M$ is finite and $\tau$ is a trace on $\M$.} Define the cone
\[
\mathscr P_\tau=\{A\in\Aff(\M):A\ge0,\ \ty(A)<\infty\}.
\]
By {Lemma~\ref{lem:formsumequals} and} Theorem~\ref{thm:additivity}, $\mathscr P_\tau$ is closed under positive sums and nonnegative scalar multiplication. {Let}
\[
\mathscr E_\tau:=\mathscr P_\tau-\mathscr P_\tau,\qquad \mathscr D_\tau:=\mathscr E_\tau+i\mathscr E_\tau.
\]

\begin{theorem}[A finite linear extension]\label{thm:linear}
The formula
\[
\widetilde\tau(A-B)=\ty(A)-\ty(B),\qquad A,B\in\mathscr P_\tau,
\]
is well-defined and real-linear on $\mathscr E_\tau$, and extends uniquely by complex linearity to $\mathscr D_\tau$. It is positive on the positive cone and agrees with $\tau$ on $\M$.
\end{theorem}
\begin{proof}
If $A-B=C-D$, then $A+D=B+C$ in {$\Aff(\M)$}. Positive additivity {(Lemma~\ref{lem:formsumequals} and Theorem~\ref{thm:additivity})} implies
\[
\ty(A)+\ty(D)=\ty(B)+\ty(C),
\]
where all terms are finite, hence the formula is independent of representation. It is additive and homogeneous. If $T=A-B\ge0$ with $A,B\in\mathscr P_\tau$, then $T+B=A$ in the affiliated algebra; positive additivity shows $\ty(T)=\ty(A)-\ty(B)\ge0$, proving positivity on this cone. Uniqueness of the decomposition into self-adjoint real and imaginary parts gives complexification. Every bounded self-adjoint $X$ is $X_+-X_-$ with both terms in $\mathscr P_\tau$, and the construction coincides with $\tau$ by Proposition~\ref{prop:elementary}.
\end{proof}

\begin{theorem}[Quadratic structure for arbitrary traces]\label{thm:L2general}
Put $\mathscr L^2_Y=\{T\in\Aff(\M):\ty(T^*T)<\infty\}$. Then $\mathscr L^2_Y$ is a complex vector space. For $S,T\in\mathscr L^2_Y$ we have $T^*S\in\mathscr D_\tau$, and
\[
\langle T,S\rangle_Y:=\widetilde\tau(T^*S)
\]
is a positive semidefinite Hermitian sesquilinear form with $\langle T,T\rangle_Y=\ty(T^*T)$. It is definite if and only if $\tau$ is faithful.
\end{theorem}
\begin{proof}
All operations are taken in the $*$-algebra $\Aff(\M)$. For closed $X$ the operator $X^*X$ is positive and self-adjoint, so it agrees with the product in $\Aff(\M)${ (Lemma~\ref{lem:posaff}(i))}.

\emph{Vector space.} In $\Aff(\M)$,
\[
(S+T)^*(S+T)+(S-T)^*(S-T)=2S^*S+2T^*T .
\]
Lemma~\ref{lem:formsumequals}, Theorem~\ref{thm:additivity} and $\ty\ge0$ give $\ty((S+T)^*(S+T))\le2\ty(S^*S)+2\ty(T^*T)$. Also $\ty((cT)^*(cT))=|c|^2\ty(T^*T)$ by Proposition~\ref{prop:elementary}.

\emph{The pairing.} The polarization identity
\[
T^*S=\tfrac14\sum_{k=0}^3 i^k\,(S+i^kT)^*(S+i^kT)
\]
holds in $\Aff(\M)$. Each term $(S+i^kT)^*(S+i^kT)$ lies in $\mathscr P_\tau$ because $\mathscr L^2_Y$ is a vector space. Hence $T^*S\in\mathscr D_\tau$, and $\widetilde\tau(T^*S)$ is defined by Theorem~\ref{thm:linear}.

\emph{Properties.} Sesquilinearity follows from linearity of $\widetilde\tau$ on $\mathscr D_\tau$ and distributivity in $\Aff(\M)$. Hermitian symmetry follows from $(T^*S)^*=S^*T$ and $\widetilde\tau(X^*)=\overline{\widetilde\tau(X)}$; the latter holds because $\widetilde\tau$ is real on $\mathscr E_\tau$. Positivity holds since $\langle T,T\rangle_Y=\widetilde\tau(T^*T)=\ty(T^*T)\ge0$.

\emph{Definiteness.} Suppose $\tau$ is faithful and $\ty(X)=0$ for some $X\ge0$. Then $\tau(\min(X,R\id))=0$ for all $R$ by Proposition~\ref{prop:clip}, so $\min(X,R\id)=0$ for all $R$, and $X=0$. If $\tau$ is not faithful, some $X\ge0$, $X\ne0$, has $\tau(X)=0$. A spectral projection $p=\mathbf 1_{[\varepsilon,\infty)}(X)\ne0$ satisfies $\varepsilon p\le X$, so $\tau(p)=0${, since $0\le\varepsilon\tau(p)\le\tau(X)=0$,} and $\langle p,p\rangle_Y=0$.
\end{proof}


\subsection{The absolute-value integrability domain}\label{sec:L1}
It is natural to introduce
\[
\mathscr L^1_Y(\M,\tau):=\{T\in\Aff(\M):\ty(|T|)<\infty\}.
\]
In an abelian algebra this is a vector space by pointwise triangle inequality and Theorem~\ref{thm:commuting}. If $\tau$ is faithful and normal, it is the conventional noncommutative $L^1(\M,\tau)$, complete for $\|T\|_1={\widetilde\tau}(|T|)$ {\cite{Nel74}}. {For an arbitrary trace we now show that $\mathscr L^1_Y$ is still a vector space, that it equals $\mathscr D_\tau$, and that it is an $\M$-bimodule on which $\widetilde\tau$ is cyclic. The triangle inequality holds up to a constant; the optimality of the constant is open.}

\begin{lemma}[An inequality for squares]\label{lem:squares}
Let $\M$ be finite and $S,T\in\Aff(\M)$. For every $t>0$,
\begin{equation}\label{eq:squares}
(S+T)^*(S+T)\le(1+t)\,S^*S+(1+t^{-1})\,T^*T\qquad\text{in }\Aff(\M).
\end{equation}
In particular, if $a,b\in\M_+$ satisfy $a^2\le S^*S$ and $b^2\le T^*T$, then $(a+b)^2\le(1+t)\,S^*S+(1+t^{-1})\,T^*T$. For $\lambda,\mu>0$, the inequality $(S+T)^*(S+T)\le\lambda S^*S+\mu T^*T$ holds for all $S$ and $T$ if and only if $\lambda^{-1}+\mu^{-1}\le1$.
\end{lemma}
\begin{proof}
In $\Aff(\M)$,
\[
(1+t)\,S^*S+(1+t^{-1})\,T^*T-(S+T)^*(S+T)=\bigl(t^{1/2}S-t^{-1/2}T\bigr)^*\bigl(t^{1/2}S-t^{-1/2}T\bigr),
\]
which lies in $\Aff(\M)_+$ by Lemma~\ref{lem:posaff}(i). Applying this to $a$ and $b$ gives $(a+b)^2\le(1+t)a^2+(1+t^{-1})b^2$, and the claim about $a$ and $b$ follows from Lemma~\ref{lem:posaff}(ii). If $\lambda^{-1}+\mu^{-1}\le1$, then $\lambda>1$; with $t=\lambda-1$ we have $\mu\ge1+t^{-1}$, so the inequality follows from~\eqref{eq:squares}. Conversely, take $S=u\id$ and $T=v\id$ with $u,v>0$. The minimum of $(\lambda u^2+\mu v^2)/(u+v)^2$ over $u,v>0$ is $(\lambda^{-1}+\mu^{-1})^{-1}$, by the Cauchy--Schwarz inequality.
\end{proof}

\begin{remark}\label{rem:squares}
In Lemma~\ref{lem:squares}, the bound for $(a+b)^2$ cannot be replaced by $(S+T)^*(S+T)$. If $T=-S$, then $(S+T)^*(S+T)=0$, although $a=b=|S|$ satisfy the hypotheses. Even for $S,T\ge0$, the hypotheses $a^2\le S^2$ and $b^2\le T^2$ do not imply $(a+b)^2\le(S+T)^2$, although they do imply $a+b\le S+T$ by Lemma~\ref{lem:heinz}. In $M_2(\mathbb C)$, take
\[
S=a=\begin{pmatrix}1&0\\0&0\end{pmatrix},\qquad T=\begin{pmatrix}1&1\\1&1\end{pmatrix},\qquad b=0;
\]
then $(S+T)^2-a^2=\left(\begin{smallmatrix}4&3\\3&2\end{smallmatrix}\right)$ has determinant $-1$. The reason is that $\lambda\mapsto\lambda^2$ is not operator monotone, whereas $\lambda\mapsto\lambda^{1/2}$ is. This is why the argument below works with squares and loses a constant.
\end{remark}

\begin{theorem}[Quasi-triangle inequality]\label{thm:L1quasi}
Let $\M$ be finite and $\tau$ a trace on $\M$. For $S,T\in\Aff(\M)$ put $\alpha=\ty(|S|)$ and $\beta=\ty(|T|)$. Then
\[
\ty(|S+T|)\le\inf_{t>0}\bigl(\sqrt{1+t}\;\alpha+\sqrt{1+t^{-1}}\;\beta\bigr)=\bigl(\alpha^{2/3}+\beta^{2/3}\bigr)^{3/2}\le\sqrt2\,(\alpha+\beta).
\]
\end{theorem}
\begin{proof}
We may assume that $\alpha,\beta<\infty$. Fix $t>0$ and let $Z_t=(1+t)|S|^2+(1+t^{-1})|T|^2\in\Aff(\M)_+$. Let $c\in\M_+$ with $c^2\le|S+T|^2$. By Lemma~\ref{lem:squares} and Lemma~\ref{lem:posaff}(ii),(iii), $c^2\le Z_t$. Corollary~\ref{cor:var}, applied to $Z_t^{1/2}$, gives $\tau(c)\le\ty(Z_t^{1/2})$. Taking the supremum over $c$ and applying Corollary~\ref{cor:var} to $S+T$, we obtain $\ty(|S+T|)\le\ty(Z_t^{1/2})$. By Lemma~\ref{lem:formsumequals}, $Z_t$ is the form sum of $(1+t)|S|^2$ and $(1+t^{-1})|T|^2$. Proposition~\ref{prop:sqrt} and homogeneity (Proposition~\ref{prop:elementary}) therefore give
\[
\ty(Z_t^{1/2})\le\ty\bigl(\sqrt{1+t}\,|S|\bigr)+\ty\bigl(\sqrt{1+t^{-1}}\,|T|\bigr)=\sqrt{1+t}\;\alpha+\sqrt{1+t^{-1}}\;\beta .
\]
If $\alpha\beta>0$, the infimum over $t$ is attained at $t=(\beta/\alpha)^{2/3}$. There $\sqrt{1+t}=(\alpha^{2/3}+\beta^{2/3})^{1/2}\alpha^{-1/3}$ and $\sqrt{1+t^{-1}}=(\alpha^{2/3}+\beta^{2/3})^{1/2}\beta^{-1/3}$, and the value is $(\alpha^{2/3}+\beta^{2/3})^{3/2}$. If $\beta=0$ (respectively $\alpha=0$), let $t\to0$ (respectively $t\to\infty$). The last inequality is the power-mean inequality $\alpha^{2/3}+\beta^{2/3}\le2^{1/3}(\alpha+\beta)^{2/3}$.
\end{proof}

\begin{corollary}[Triangle inequality on average]\label{cor:average}
Let $\M$ be finite and $\tau$ a trace on $\M$. For $S,T\in\Aff(\M)$,
\[
\tfrac12\bigl(\ty(|S+T|)+\ty(|S-T|)\bigr)\le\ty\bigl((S^*S+T^*T)^{1/2}\bigr)\le\ty(|S|)+\ty(|T|).
\]
\end{corollary}
\begin{proof}
Put $A=\frac12(S+T)$ and $B=\frac12(S-T)$, so that $2A^*A+2B^*B=S^*S+T^*T$. If $a,b\in\M_+$ satisfy $a^2\le A^*A$ and $b^2\le B^*B$, Lemma~\ref{lem:squares} with $t=1$ gives $(a+b)^2\le S^*S+T^*T$. By Corollary~\ref{cor:var}, $\tau(a)+\tau(b)=\tau(a+b)\le\ty\bigl((S^*S+T^*T)^{1/2}\bigr)$. Taking suprema over $a$ and $b$, Corollary~\ref{cor:var} and homogeneity give the first inequality, since $\ty(|A|)=\frac12\ty(|S+T|)$ and $\ty(|B|)=\frac12\ty(|S-T|)$. The second inequality is Proposition~\ref{prop:sqrt} with $X=S^*S$ and $Y=T^*T$, together with Lemma~\ref{lem:formsumequals}.
\end{proof}

\begin{theorem}[Integrable operators]\label{thm:L1}
Let $\M$ be finite and $\tau$ a trace on $\M$.
\begin{enumerate}[label=\textup{(\alph*)},leftmargin=*]
\item $\mathscr L^1_Y$ is a complex vector space, $\mathscr L^1_Y=\mathscr D_\tau$, and $d(S,T)=\ty(|S-T|)^{2/3}$ is a translation-invariant pseudometric on $\mathscr L^1_Y$. It is a metric iff $\tau$ is faithful.
\item For $T\in\Aff(\M)$ and $X\in\M$, 
one has $|T^*|=U|T|U^*$ with $U\in\M$ unitary, $\ty(|T^*|)=\ty(|T|)$, 
$\ty(|XT|)\le\|X\|\,\ty(|T|)$ and $\ty(|TX|)\le\|X\|\,\ty(|T|)$. In particular, $\mathscr L^1_Y$ is an $\M$-bimodule.
\item $|\widetilde\tau(T)|\le\sqrt2\,\ty(|T|)$ for $T\in\mathscr L^1_Y$.
\item $\widetilde\tau(XT)=\widetilde\tau(TX)$ for $X\in\M$ and $T\in\mathscr L^1_Y$.
\end{enumerate}
\end{theorem}
\begin{proof}
Since $\M$ is finite, every $T\in\Aff(\M)$ has a polar decomposition $T=U|T|$ with $U\in\M$ unitary~\cite[Theorem~2.12]{Liu11}. Then $TT^*=U|T|^2U^*$, so $|T^*|=U|T|U^*$ and $\ty(|T^*|)=\ty(|T|)$ by Proposition~\ref{prop:elementary}.

(a) By Theorem~\ref{thm:L1quasi} and $|cT|=|c|\,|T|$, $\mathscr L^1_Y$ is a vector space. The function $N(T)=\ty(|T|)^{2/3}$ satisfies $N(S+T)\le N(S)+N(T)$ and $N(-T)=N(T)$, so $d$ is a translation-invariant pseudometric. If $T\in\mathscr L^1_Y$, then $T^*\in\mathscr L^1_Y$ by the preceding paragraph, so $\operatorname{Re}T,\operatorname{Im}T\in\mathscr L^1_Y$. Their positive and negative parts are dominated by their absolute values, so they lie in $\mathscr P_\tau$. Hence $T\in\mathscr D_\tau$. Conversely, $\mathscr D_\tau$ is spanned by $\mathscr P_\tau\subseteq\mathscr L^1_Y$.

If $\tau$ is faithful and $\ty(|R|)=0$, then $\tau(\min(|R|,\id))=0$ by Proposition~\ref{prop:clip}, so $\min(|R|,\id)=0$ and $R=0$; conversely, if $d$ is a metric and $X\in\M_+$ satisfies $\tau(X)=0$, then $d(X,0)=\tau(X)^{2/3}=0$, so $X=0$; hence $\tau$ is faithful.

(b) The first assertion was shown above. Let $D=(\|X\|^2\id-X^*X)^{1/2}\in\M_+$. In $\Aff(\M)$,
\[
\|X\|^2\,T^*T-(XT)^*(XT)=T^*D^2T=(DT)^*(DT),
\]
so $|XT|^2\le\bigl(\|X\|\,|T|\bigr)^2$ by Lemma~\ref{lem:posaff}(i). If $c\in\M_+$ and $c^2\le|XT|^2$, then $c^2\le\bigl(\|X\|\,|T|\bigr)^2$ by Lemma~\ref{lem:posaff}(ii),(iii), and Corollary~\ref{cor:var} gives $\tau(c)\le\ty(\|X\|\,|T|)=\|X\|\,\ty(|T|)$. Taking the supremum over $c$, Corollary~\ref{cor:var} yields $\ty(|XT|)\le\|X\|\,\ty(|T|)$. Finally, $\ty(|TX|)=\ty(|X^*T^*|)\le\|X\|\,\ty(|T^*|)=\|X\|\,\ty(|T|)$.

(c) Recall from the proof of Theorem~\ref{thm:L2general} that $\widetilde\tau(R^*)=\overline{\widetilde\tau(R)}$ for $R\in\mathscr D_\tau$. Let $T\in\mathscr L^1_Y$, choose $\theta\in\mathbb R$ with $e^{i\theta}\widetilde\tau(T)=|\widetilde\tau(T)|$, and put $H=\operatorname{Re}(e^{i\theta}T)\in\mathscr L^1_Y$. Then $H=H_+-H_-$ with $H_\pm\in\mathscr P_\tau$ and $H_+\le|H|$, so
\[
|\widetilde\tau(T)|=\operatorname{Re}\widetilde\tau(e^{i\theta}T)=\widetilde\tau(H)=\ty(H_+)-\ty(H_-)\le\ty(|H|).
\]
Since $\ty(|e^{i\theta}T|)=\ty(|e^{-i\theta}T^*|)=\ty(|T|)$ by (b), Theorem~\ref{thm:L1quasi} gives
\[
\ty(|H|)=\tfrac12\,\ty\bigl(|e^{i\theta}T+e^{-i\theta}T^*|\bigr)\le\tfrac12\bigl(2\,\ty(|T|)^{2/3}\bigr)^{3/2}=\sqrt2\,\ty(|T|).
\]

(d) Let $A\in\mathscr P_\tau$ and $A_R=\min(A,R\id)\in\M$. By Theorem~\ref{thm:commuting} and Proposition~\ref{prop:clip}, $\ty(A-A_R)=\ty(A)-\tau(A_R)\to0$. By (b) and (c),
\[
|\widetilde\tau(XA)-\tau(XA_R)|=|\widetilde\tau(X(A-A_R))|\le\sqrt2\,\|X\|\,\ty(A-A_R)\to0,
\]
and likewise for $AX$. Since $\tau(XA_R)=\tau(A_RX)$, we get $\widetilde\tau(XA)=\widetilde\tau(AX)$. The claim follows by linearity, since $\mathscr L^1_Y=\mathscr D_\tau$ is spanned by $\mathscr P_\tau$.
\end{proof}
\section{Faithful normal traces}\label{sec:normal}
{In this section} the \emph{standing assumptions} are that $\M\subseteq\mathcal B(H)$ is an arbitrary von Neumann algebra and $\tau:\M\to\mathbb C$ is a faithful normal tracial state, $\tau(\id)=1$. We do \emph{not} assume finiteness of $\M$: it follows from faithfulness.

\begin{proposition}[A faithful tracial state forces finiteness]\label{prop:automaticfinite}
Let $\M$ be any von Neumann algebra admitting a positive faithful tracial functional $\tau$ with $0<\tau(\id)<\infty$. Then $\M$ is finite in the Murray--von Neumann sense: every isometry in $\M$ is unitary. Normality of $\tau$ is not required for this conclusion.
\end{proposition}
\begin{proof}
If $v\in\M$ satisfies $v^*v=\id$, then $p=\id-vv^*$ is a projection. Traciality gives
\[
\tau(p)=\tau(\id)-\tau(vv^*)=\tau(\id)-\tau(v^*v)=0.
\]
Faithfulness implies $p=0$, so $vv^*=\id$. Thus $v$ is unitary.
\end{proof}

For positive affiliated $A$, we have (Lemma~\ref{lem:normal-spectral})
\[
\widetilde\tau(A)=\int_{[0,\infty)}\lambda\,d\tau(E_A(\lambda)) =\lim_{h\downarrow0}\tau\bigl(A(\id+hA)^{-1}\bigr).
\]
 We have $\mathscr D_\tau=\mathscr L^1_Y=L^1(\M,\tau)$. 
More generally, let $1\leq p<\infty$, let $\mathscr L^p_Y$ be the set of $T\in\Aff(\M)$ such that $\|T\|_{p,\tau}:=\tau_Y(|T|^p)^{1/p}<\infty$.
\begin{theorem}[\cite{FK86}]\label{thm:L2normal}
{Let $\tau$ be a faithful normal tracial state on $\M$. Then $\mathscr L^p_Y=L^p(\M,\tau)$ with equal norms; these spaces are complete. 
}
\end{theorem}

 For $A\in L^1$ and $D_h\in\M$,
\[
|\tau(D_h)-\widetilde\tau(A)|\leq\|D_h-A\|_{1,\tau};
\]
{this estimate allows  to recover the trace of the integrable unbounded operators via Chernoff approximation without an operator-norm limit.}

\begin{lemma}[The $L^1$--$L^2$ inequality for a normalized normal trace]\label{lem:L1L2}
For $T\in L^2(\M,\tau)$ with $\tau(\id)=1$, one has $\|T\|_{1,\tau}\le \|T\|_{2,\tau}$.
\end{lemma}
\begin{proof}
Apply Cauchy--Schwarz for the normal trace to $\id$ and $|T|$, first with $|T|$ replaced by its bounded clipping $\min(|T|,R\id)$:
\[
\tau(\min(|T|,R\id))\le\tau(\id)^{1/2}\tau(\min(|T|,R\id)^2)^{1/2}\le{\widetilde\tau}(|T|^2)^{1/2}.
\]
Take $R\to\infty$ by normal monotone convergence.
\end{proof}

\begin{lemma}[Bounded strong convergence in finite tracial $L^p$]\label{lem:stronglp}
Suppose $X_n,X\in\M$, $X_n\to X$ strongly and $\sup_n\|X_n\|<\infty$. Then $\|X_n-X\|_{p,\tau}\to0$ for each $1\leq p<\infty$.
\end{lemma}
\begin{proof}
Set $Y_n=X_n-X$. Since $Y_n\to0$ strongly and $\sup_n\|Y_n\|\leq C$, also $Y_n^*Y_n\to0$ strongly: $\|Y_n^*Y_n\xi\|\leq C\|Y_n\xi\|$. {Hence $Y_n^*Y_n\to0$ weakly, with $\|Y_n^*Y_n\|\le C^2$. A normal positive functional is $\sigma$-weakly continuous, and the $\sigma$-weak and weak operator topologies agree on bounded sets;\footnote{See~\cite[Chap.~7]{KR86}, in particular~\cite[Theorem~7.1.12]{KR86} for further details.} therefore $\tau(Y_n^*Y_n)\to0$.} For $p\leq2$, finite-trace H\"older gives $\|Y_n\|_p\leq\|Y_n\|_2$. For $p>2$, $|Y_n|^p\leq C^{p-2}|Y_n|^2$, proving the claim.
\end{proof}

\begin{corollary}[Iterated Chernoff recovery of the extended trace]\label{cor:tracechernoff}
Let $A\geq0$ be self-adjoint affiliated with $\M$. Suppose $F(h)\in\M$ and, for each fixed $t>0$, $F(t/n)^n\to e^{-tA}$ strongly as $n\to\infty$. Then for $1\leq p<\infty$,
\[
\|F(t/n)^n-e^{-tA}\|_{p,\tau}\to0,
\]
and
\begin{equation}\label{eq:iterated}
\widetilde\tau(A)=\lim_{t\downarrow0}\lim_{n\to\infty}\frac1t\tau\bigl(\id-F(t/n)^n\bigr)\in[0,\infty].
\end{equation}
\end{corollary}
\begin{proof}
For fixed $t$, strong convergence and the uniform boundedness principle give a uniform norm bound on $F(t/n)^n$. Apply Lemma~\ref{lem:stronglp}. The inner limit of the scalar expression is $t^{-1}\tau(\id-e^{-tA})$. Its limit as $t\downarrow0$ is $\widetilde\tau(A)$ by {Proposition~\ref{prop:exp}, or by} spectral monotone convergence, since $(1-e^{-t\lambda})/t\uparrow\lambda$.
\end{proof}

\begin{remark}[Order of limits]
Equation~\eqref{eq:iterated} is an \emph{iterated}-limit identity. It supplies no uniform rate as $t\downarrow0$ and does not justify an arbitrary diagonal choice $n=n(t)$.
Observe also that neither positivity nor self-adjointness of the individual $F(h)$ is required for the iterated identity. Moreover, faithfulness is not used in Corollary~\ref{cor:tracechernoff}.
\end{remark}


\appendix

\section{Affiliated operators: proofs}\label{app:aff}
\subsection{Proof of Proposition~\ref{prop:spec-aff}}\label{app:specaff}
$\Longrightarrow)$. Fix a unitary $U\in\M'$ and a Borel $B\subseteq\mathbb R$. Then, $UAU^*=A$. {The map $B\mapsto UE_A(B)U^*$ is a projection-valued measure on $\mathbb R$, and  $UAU^*=\int_{\sigma(A)}\lambda\,d(UE_AU^*)(\lambda)$. By uniqueness in the spectral theorem~\cite[Ch.~VIII]{RS80}, $E_A(B)=E_{UAU^*}(B)=UE_A(B)U^*$.}
%
The unitaries of the von Neumann algebra $\M'$ linearly span $\M'$; thus $E_A(B)$ commutes with \emph{every} element of $\M'$. The bicommutant theorem implies $E_A(B)\in\M''=\M$.

\emph{Bounded Borel calculus.} For a bounded simple function $f=\sum_{j=1}^N c_j\mathbf1_{B_j}$, we have $f(A)=\sum_jc_jE_A(B_j)\in\M$. Every bounded Borel function is a uniform limit of bounded simple functions, and the Borel calculus is contractive: $\|f_n(A)-f(A)\|\leq\|f_n-f\|_\infty$. Sequential closedness of $\M$ therefore yields $f(A)\in\M$.

$\Longleftarrow)$. Assume $E_A(B)\in\M$ for every Borel $B$. Every unitary $U\in\M'$ then commutes with these projections. For $R>0$ let\footnote{$A_R\in \M$ by the preceding paragraph applied to $f(\lambda)=\lambda\mathbf 1_{[-R,R]}(\lambda)$.}
\[
A_R=\int_{[-R,R]}\lambda\,dE_A(\lambda)=A E_A([-R,R])\in\M,
\]
 In particular, $UA_R=A_RU$. The spectral theorem yields the following characterization:
\[
\Dom(A)=\Bigl\{x\in H:\int_{\mathbb R}\lambda^2\,d\|E_A(\lambda)x\|^2<\infty\Bigr\}=\Bigl\{x:\sup_{R>0}\|A_Rx\|<\infty\Bigr\}.
\]
For $x\in\Dom(A)$, it follows $A_Rx\to Ax$. Therefore $A_RUx=UA_Rx\to UAx$ as $R\to\infty$ and $\sup_R\|A_RUx\|\leq\|Ax\|$. The above domain characterization shows $Ux\in\Dom(A)$ and $AUx=UAx$. Repeating this with $U^*$ gives $U\Dom(A)=\Dom(A)$. Hence $A\eta\M$.

{\emph{Unbounded Borel calculus.} Let $f:\mathbb R\to\mathbb C$ be Borel. Then $f(A)=\int f\,dE_A$ is closed, densely defined and normal~\cite[Ch.~VIII]{RS80}. For a unitary $U\in\M'$, the first paragraph gives $UE_A(B)U^*=E_A(B)$ for every Borel $B$, so $Uf(A)U^*=\int f\,d(UE_AU^*)=f(A)$, including domains. Hence $f(A)\eta\M$. If $f$ is real-valued, $f(A)$ is self-adjoint, and if $f\ge0$, then $\langle f(A)x,x\rangle=\int f\,d\langle E_Ax,x\rangle\ge0$ for $x\in\Dom f(A)$.} \qed

\subsection{Proof of Lemma~\ref{lem:formsumequals}}\label{app:formsum}
By the Murray--von Neumann theorem on the affiliated $*$-algebra of a finite von Neumann algebra~\cite[Propositions~3.7(i) and~3.10]{Liu11}, $\Dom A\cap\Dom B$ is dense, $A+B$ is closable there, and ${(A\splus B)^*=A^*\splus B^*=A\splus B}$. Hence $\mathcal D$ is dense. The form $t[y]=\|A^{1/2}y\|^2+\|B^{1/2}y\|^2$ on $\mathcal D$ is closed, being a sum of closed nonnegative forms and defines a self-adjoint operator $S\geq 0$. We have $S=A\splus B\ge0$ with $\Dom(S^{1/2})=\mathcal D$ and $\|S^{1/2}y\|^2=t[y]$. Indeed, if $x\in\Dom A\cap\Dom B$ and $y\in\mathcal D$, then $t(y,x)=\langle y,(A+B)x\rangle$. By the first representation theorem, $x\in\Dom S$ and $Sx=(A+B)x$. Thus ${A\splus B}\subseteq S$, so $S=S^*\subseteq {(A\splus B)^*=A\splus B}$, and equality holds.

Affiliation: for a unitary $u\in\M'$ one has $uA^{1/2}\subseteq A^{1/2}u$ and $uB^{1/2}\subseteq B^{1/2}u$. Hence $u\mathcal D=\mathcal D$ and $t[uy]=t[y]$. The operator associated with a closed form is determined by the form, so $uSu^*=S$. \qed

\subsection{Two affiliation criteria}\label{app:affcrit}
\begin{lemma}[Two affiliation criteria]\label{lem:affcrit}
Let $\M$ be any von Neumann algebra.
\begin{enumerate}[label=\textup{(\roman*)},leftmargin=*]
\item Let $T$ be closed and densely defined, with a core $\mathcal D_0$. Suppose that, for every unitary $U\in\M'$ and every $x\in\mathcal D_0$, we have $Ux\in\Dom(T)$ and $TUx=UTx$. Then $T\eta\M$.
\item Let $T$ be closed and affiliated with $\M$, and let $X\in\M$ with $\Ran X\subseteq\Dom(T)$. Then $TX$ is everywhere defined and $TX\in\M$.
\end{enumerate}
\end{lemma}
\begin{proof}
(i) Let $y\in\Dom(T)$. Choose $y_n\in\mathcal D_0$ with $y_n\to y$ and $Ty_n\to Ty$. Then $Uy_n\to Uy$ and $TUy_n=UTy_n\to UTy$. Since $T$ is closed, $Uy\in\Dom(T)$ and $TUy=UTy$. Thus $U\Dom(T)\subseteq\Dom(T)$, and applying this to $U^*$ gives equality.

(ii) $TX$ is everywhere defined because $\Ran X\subseteq\Dom(T)$. It is closed: if $x_n\to x$ and $TXx_n\to y$, then $Xx_n\to Xx$, and closedness of $T$ gives $Xx\in\Dom(T)$ and $TXx=y$~\cite[Lemma~4.10]{Liu11}. By the closed graph theorem, $TX$ is bounded. Let $U\in\M'$ be unitary and $x\in H$. Then $UXx=XUx$, and, since $T\eta\M$ and $Xx\in\Dom(T)$, also $TUXx=UTXx$. Hence $(TX)Ux=T(XUx)=T(UXx)=U(TX)x$. The unitaries of $\M'$ span $\M'$, so $TX\in\M''=\M$.
\end{proof}

\subsection{Proof of Lemma~\ref{lem:posaff}}\label{app:posaff}
\emph{Self-adjointness.} For $x,y\in\Dom Z$, the sesquilinear forms $\langle Zx,y\rangle$ and $\langle x,Zy\rangle$ have the same quadratic form, because $\langle Zx,x\rangle$ is real. By polarization they coincide, so $Z\subseteq Z^*$. A symmetric operator affiliated with a finite von Neumann algebra is self-adjoint~\cite[Proposition~3.5]{Liu11}. If $\lambda<0$, then $\|(Z-\lambda)x\|\,\|x\|\ge\langle(Z-\lambda)x,x\rangle\ge|\lambda|\,\|x\|^2$ for $x\in\Dom Z$. Thus $Z-\lambda$ is injective with closed range, and its range is dense because $\ker(Z^*-\lambda)=\ker(Z-\lambda)=\{0\}$. Hence $\lambda\notin\sigma(Z)$.

(i) By von Neumann's theorem~\cite[Theorem~2.4(v)]{Liu11}, $T^*T$ is self-adjoint, and $\langle T^*Tx,x\rangle=\|Tx\|^2\ge0$. It is affiliated with $\M$, because $U^*T^*TU=(U^*TU)^*(U^*TU)=T^*T$ for every unitary $U\in\M'$. Since $T^*T$ is closed, it equals its closure $T^*\sdot T$.

(ii) Let $Z_1,Z_2\in\Aff(\M)_+$. On $\Dom Z_1\cap\Dom Z_2$, which is a core for $Z_1\splus Z_2$, we have $\langle(Z_1+Z_2)x,x\rangle\ge0$, and this inequality passes to the closure. Closure under nonnegative multiples is clear. If $W\le V\le Z$, then $Z-W=(Z-V)+(V-W)\in\Aff(\M)_+$ by the algebraic rules of $\Aff(\M)$~\cite[Theorem~3.11]{Liu11}.

(iii) If $W,Z\in\Aff(\M)_+$ and $D:=Z-W\in\Aff(\M)_+$, then $Z=W\splus D$, and $Z$ is the form sum of $W$ and $D$ by Lemma~\ref{lem:formsumequals}. Hence $\Dom(Z^{1/2})\subseteq\Dom(W^{1/2})$ and $\|Z^{1/2}x\|^2=\|W^{1/2}x\|^2+\|D^{1/2}x\|^2\ge\|W^{1/2}x\|^2$ for $x\in\Dom(Z^{1/2})$. Conversely, let $W\in\M_+$ with $W\le Z$ as quadratic forms. The operator $Z-W$, with domain $\Dom Z$, is closed and affiliated with $\M$, so it is the difference in $\Aff(\M)$; and $\langle(Z-W)x,x\rangle\ge0$ for $x\in\Dom Z$. \qed

\section{Proofs for the square-root inequality}\label{app:sqrt}

Throughout this appendix $0<\varepsilon<M<\infty$, and
\[
g_{\varepsilon,M}(\lambda)=\frac1\pi\int_\varepsilon^M s^{-1/2}\frac{\lambda}{\lambda+s}\,ds=\frac1\pi\int_\varepsilon^M s^{-3/2}f_{1/s}(\lambda)\,ds\qquad(\lambda\ge0).
\]

\begin{lemma}[Truncated square root]\label{lem:gtrunc}
\begin{enumerate}[label=\textup{(\alph*)},leftmargin=*]
\item $g_{\varepsilon,M}$ is continuous and increasing, and $0\le g_{\varepsilon,M}(\lambda)\le\min\bigl(\lambda^{1/2},2M^{1/2}/\pi\bigr)$.
\item For $\lambda\ge0$, $\lambda^{1/2}-g_{\varepsilon,M}(\lambda)\le\pi^{-1}\bigl(2\varepsilon^{1/2}+2\lambda M^{-1/2}\bigr)$. For $R>0$,
\[
\min(\lambda^{1/2},R)\le g_{\varepsilon,M}(\lambda)+\delta,\qquad \delta=\pi^{-1}\bigl(2\varepsilon^{1/2}+2R^2M^{-1/2}\bigr).
\]
\item Let $Z\ge0$ be self-adjoint and affiliated with $\M$. Then $s\mapsto s^{-3/2}f_{1/s}(Z)$ is norm continuous from $[\varepsilon,M]$ to $\M$, its Riemann sums converge in norm, and
\begin{gather*}
g_{\varepsilon,M}(Z)=\frac1\pi\int_\varepsilon^M s^{-3/2}f_{1/s}(Z)\,ds,\\
\tau\bigl(g_{\varepsilon,M}(Z)\bigr)=\frac1\pi\int_\varepsilon^M s^{-3/2}\,\tau\bigl(f_{1/s}(Z)\bigr)\,ds
\end{gather*}
for every positive linear functional $\tau$ on $\M$.
\end{enumerate}
\end{lemma}
\begin{proof}
(a) The integrand increases with $\lambda$ and lies between $0$ and $s^{-1/2}$, which gives monotonicity, continuity and the bound $2M^{1/2}/\pi$. The substitution $s=\lambda u$ and $\int_0^\infty u^{-1/2}(1+u)^{-1}\,du=\pi$ give $\pi\lambda^{1/2}=\int_0^\infty s^{-1/2}\frac{\lambda}{\lambda+s}\,ds$, hence $g_{\varepsilon,M}(\lambda)\le\lambda^{1/2}$.

(b) By (a),
\begin{align*}
\lambda^{1/2}-g_{\varepsilon,M}(\lambda)&=\frac1\pi\Bigl(\int_0^\varepsilon+\int_M^\infty\Bigr)s^{-1/2}\frac{\lambda}{\lambda+s}\,ds\\
&\le\frac1\pi\Bigl(\int_0^\varepsilon s^{-1/2}\,ds+\lambda\int_M^\infty s^{-3/2}\,ds\Bigr),
\end{align*}
which is the first claim. If $\lambda\le R^2$, it gives $\min(\lambda^{1/2},R)=\lambda^{1/2}\le g_{\varepsilon,M}(\lambda)+\delta$. If $\lambda>R^2$, monotonicity gives $g_{\varepsilon,M}(\lambda)\ge g_{\varepsilon,M}(R^2)\ge R-\delta$.

(c) Since $f_{1/s}(\lambda)=s\lambda/(s+\lambda)$, for $s,s'>0$ and $\lambda\ge0$ we have
\[
\bigl|f_{1/s}(\lambda)-f_{1/s'}(\lambda)\bigr|=\frac{\lambda^2\,|s-s'|}{(s+\lambda)(s'+\lambda)}\le|s-s'| .
\]
Moreover $0\le f_{1/s}\le s\le M$, so the functions $\varphi_s=s^{-3/2}f_{1/s}$, $s\in[\varepsilon,M]$, satisfy $\sup_{\lambda\ge0}|\varphi_s(\lambda)-\varphi_{s'}(\lambda)|\le\omega(|s-s'|)$ for some $\omega$ with $\omega(r)\to0$ as $r\to0$. By Proposition~\ref{prop:spec-aff}, $\varphi_s(Z)\in\M$, and $\|\psi(Z)\|\le\sup_{\lambda\ge0}|\psi(\lambda)|$ for bounded Borel $\psi$; hence $s\mapsto\varphi_s(Z)$ is norm continuous whence is Riemann-integrable. Positive functionals are norm continuous (Definition~\ref{def:trace}), so $\tau$ commutes with the integral.
\end{proof}

\begin{proof}[Proof of Proposition~\ref{prop:sqrt}]
Put $Z=X\splus Y$. By Lemma~\ref{lem:gtrunc}(c) and Theorem~\ref{thm:additivity}(b) with $h=1/s$,
\[
\begin{aligned}
\tau\bigl(g_{\varepsilon,M}(Z)\bigr)&=\frac1\pi\int_\varepsilon^Ms^{-3/2}\,\tau\bigl(f_{1/s}(Z)\bigr)\,ds\\
&\le\frac1\pi\int_\varepsilon^Ms^{-3/2}\bigl[\tau\bigl(f_{1/s}(X)\bigr)+\tau\bigl(f_{1/s}(Y)\bigr)\bigr]\,ds\\
&=\tau\bigl(g_{\varepsilon,M}(X)\bigr)+\tau\bigl(g_{\varepsilon,M}(Y)\bigr).
\end{aligned}
\]
By Lemma~\ref{lem:gtrunc}(a), $g_{\varepsilon,M}(X)\in\M_+$ and $g_{\varepsilon,M}(X)\le X^{1/2}$, so $\tau(g_{\varepsilon,M}(X))\le\ty(X^{1/2})$ by~\eqref{eq:inner}; likewise for $Y$. By Lemma~\ref{lem:gtrunc}(b) and the functional calculus, $\min(Z^{1/2},R\id)\le g_{\varepsilon,M}(Z)+\delta\id$. Hence
\[
\tau\bigl(\min(Z^{1/2},R\id)\bigr)\le\ty(X^{1/2})+\ty(Y^{1/2})+\delta\,\tau(\id).
\]
Let $\varepsilon\to0$ and $M\to\infty$, so that $\delta\to0$. Taking the supremum over $R$ and using Proposition~\ref{prop:clip} completes the proof.\footnote{Observe that only norm limits have been taken inside $\tau$.}
\end{proof}

\begin{lemma}[Heinz inequality for a bounded minorant]\label{lem:heinz}
Let $c\in\mathcal B(H)$ with $c\ge0$, and let $Z\ge0$ be self-adjoint with $c^2\le Z$. Then $c\le Z^{1/2}$.
\end{lemma}
\begin{proof}
For $s>0$, $c^2+s\id\le Z+s\id$ gives $(Z+s\id)^{-1}\le(c^2+s\id)^{-1}$, that is, $c^2(c^2+s\id)^{-1}\le Z(Z+s\id)^{-1}$. Multiply by $\pi^{-1}s^{-1/2}$ and integrate over $[\varepsilon,M]$; by Lemma~\ref{lem:gtrunc}(c) with $\M=\mathcal B(H)$, this gives $g_{\varepsilon,M}(c^2)\le g_{\varepsilon,M}(Z)\le Z^{1/2}$ in the sense of quadratic forms. By Lemma~\ref{lem:gtrunc}(b), $\|c-g_{\varepsilon,M}(c^2)\|\le\pi^{-1}\bigl(2\varepsilon^{1/2}+2\|c\|^2M^{-1/2}\bigr)\to0$ as $\varepsilon\to0$ and $M\to\infty$. Hence $\langle cx,x\rangle\le\|Z^{1/4}x\|^2$ for every $x\in\Dom(Z^{1/4})$.
\end{proof}

\section*{Acknowledgments}
The author thanks Prof. Wilfrid Gangbo for a fruitful initial discussion that shaped this project.

\section*{AI Disclosure.} Claude was used for mathematical discussions and editorial assistance. The author assumes full responsibility for the content.

\bibliographystyle{plain}
\bibliography{references}

\end{document}